\documentclass{amsart}

\usepackage{graphicx}
\usepackage{amssymb,amscd, amsthm, latexsym, mathrsfs, epsfig}

\usepackage{xcolor}
\usepackage{todonotes}

\newcommand\PP{\mathbb P}

\newcommand\C{\mathbb C}

\newcommand\R{\mathbb R}
\newcommand\Z{\mathbb Z}

\newcommand{\TT}{\mathbb{T}}

\newcommand{\calH}{\mathcal{H}}

\newcommand{\M}{\mathcal{M}}

\newcommand{\cS}{\mathcal{S}}

\newcommand{\cP}{\mathcal{P}}

\newcommand{\cW}{ \mathcal{W}}
\newcommand{\cL}{ \mathcal{L}}
\newcommand{\cJ}{ \mathcal{J}}

\newcommand{\fD}{\mathfrak{D}}
\newcommand{\fc}{\mathfrak{c}}

\newcommand{\fh}{\mathfrak{h}}
\newcommand{\fe}{\mathfrak{e}}
\newcommand{\ft}{\mathfrak{t}}
\newcommand{\fz}{\mathfrak{z}}

\newcommand\Aut{\operatorname{Aut}}

\newcommand\Lie{\operatorname{Lie}}
\newcommand\Asymp{\operatorname{Asymp}}
\newcommand\ii{\operatorname{i}}

\newcommand\bra{\langle}
\newcommand\ket{\rangle}

\newcommand\res{\operatorname{res}}
\newcommand\Res{\operatorname{Res}}
\newcommand\IR{\operatorname{IR}}

\newcommand\Int{\operatorname{Int}}

\newcommand\opw{\operatorname{w}}
\newcommand\opL{\operatorname{L}}
\newcommand\opH{\operatorname{H}}
\newcommand\LT{\operatorname{LT}}
\newcommand\LR{\operatorname{LR}}
\newcommand\mr{\operatorname{MR}}
\newcommand\mt{\operatorname{MT}}

\newcommand{\del}{\partial}

\newcommand{\jk}{\operatorname{JK}}

\makeatletter \@addtoreset{equation}{section} \makeatother

\newtheorem{thm}{Theorem}[section]
\newtheorem{prop}[thm]{Proposition}
\newtheorem{lem}[thm]{Lemma}
\newtheorem{cor}[thm]{Corollary}

\theoremstyle{definition}
\newtheorem{definition}[thm]{Definition}
\newtheorem{exm}[thm]{Example}
\newtheorem{rmk}[thm]{Remark}

\newcommand\narrowdots{\hbox to 1em{.\hss.\hss.}}

\title{Scattering diagrams and Jeffrey-Kirwan residues}

\author{Sara Angela Filippini}
\email{saraangela.filippini@unisalento.it}
\address{Universit\`a del Salento, Dipartimento di Matematica e Fisica ``Ennio De Giorgi", Lecce, Italy} 

\author{Jacopo Stoppa}
\email{jstoppa@sissa.it}
\address{SISSA, via Bonomea 265, 34136 Trieste, Italy} 
\address{Institute for Geometry and Physics (IGAP), Via Beirut 2-4, 34014 Trieste, Italy}
\begin{document}
	
	\begin{abstract} We show that the consistent completion of an initial scattering diagram in $M_{\R}$ (for a finite rank lattice $M$) can be expressed quite generally in terms of the Jeffrey-Kirwan residues of certain explicit meromorphic forms, by using the Maurer-Cartan asymptotic analysis developed by Chan-Leung-Ma and Leung-Ma-Young. A similar result holds for the associated theta functions.
	\end{abstract}
	
	\maketitle
	
\section{Introduction}

Residues of meromorphic forms associated with hyperplane arrangements, in the sense of Jeffrey-Kirwan \cite{jk} (as studied by Brion, Szenes and Vergne \cite{brionvergne, szenesvergne}, denoted by $\jk$ in the following), appear naturally in computations in physics (see e.g. \cite{beau, beht, cordova}) and enumerative geometry (see e.g. \cite{KimUeda, ontani, szenesvergne}, which also contain more complete lists of references). We will review this notion briefly, in a special case, in Section \ref{JKSec} of the present paper (see also \cite{ontani} for a recent detailed exposition). Here we give two motivating examples.
\begin{exm}[\cite{KimUeda}, Section 2]\label{GLSMExm} Let $W$ denote a finite dimensional complex vector space endowed with commuting actions of $G\times\C^*$ and $H$, where $G$ is a reductive algebraic group of rank $r$ and $H$ is an algebraic torus. Fix a maximal torus $T \subset G$, and write $(\rho_i, r_i, \nu_i) \in \ft^{\vee} \oplus \Z \oplus \fh^{\vee}$ for the $T\times\C^*\times H$-weight of $W_i$, where $W = \bigoplus^N_{i=1} W_i$ is the weight space decomposition with respect to $T\times\C^*$. Write $t'$ for an element of $\fz\otimes\C$, where $\fz = \Lie(Z(G))$. Introduce the rational functions of $x \in \ft$, for fixed $\lambda \in \fh$, given by
\begin{align}\label{QuasiMapZ}
\nonumber& Z^{\operatorname{vec}}_d(x) = \prod_{\alpha \in \Delta_+}(-1)^{\alpha(d)+1}\alpha^2(x),\\
\nonumber & Z^{\operatorname{mat}}_d(x) = \prod^N_{i=1} (\rho_i(x) + \nu_i(\lambda))^{r_i-\rho_i(d)-1},\\
& Z_d(x) = Z^{\operatorname{vec}}_d(x)Z^{\operatorname{mat}}_d(x),
\end{align}
{where $\Delta_+$ denotes the set of positive roots.} Fix a polynomial $P(x)$ which is invariant under the Weyl group $\cW = N(T)/T$. Write $\operatorname{P}^{\vee}$ for the coweight lattice of $G$ and $\fc \subset \fz^{\vee}$ for the K\"ahler cone of the GIT quotient $W/\!\!/_{t'}G$. Then, according to \cite{KimUeda}, Conjecture 10.10, under suitable assumptions, for $d \in \operatorname{P}^{\vee}$, {we have
\begin{equation}\label{KimConj}
\bra P \ket_{W/\!\!/_{t'}G, d} = \frac{1}{|\cW|} \jk_{\fc}(Z_d(x)P(x)),
\end{equation}
where the left hand side denotes the \emph{quasimap invariant} of $W/\!\!/_{t'}G$ with insertion $P(x)$ and degree $d$ (where $t'$ is thought of as a K\"ahler parameter), that is, an enumerative invariant for rational curves in $W/\!\!/_{t'}G$ discussed in \cite{KimUeda}, Section 10, while the left hand side is the Jeffrey-Kirwan residue of the meromorphic function $Z_d(x)P(x)$ along the hyperplane arrangement given by the singular locus of $Z_d(x)$ (by definition, this also depends on the choice of the stability parameter $t'$).} 

In the case when $G$ is itself a torus, this was a conjecture of Batyrev and Materov (\cite{BatyrevMaterov}, {Conjecture 4.6}), proved by Szenes and Vergne (\cite{szenesvergne}, {Theorem 4.1}).
\end{exm}
\begin{exm}[\cite{beau}, Section 3]\label{QuiverExm}	Fix a quiver $Q$ {with vertices $Q_0$ and arrows $Q_1$}, without loops or oriented cycles, together with a dimension vector $d$ {(see e.g. \cite{MeinhardtReineke}, Section 2.1 for some background on quivers and their representations)}. The gauge group is a product of unitary groups $\prod_{i \in Q_0} \operatorname{U}(d_i)$, modulo the diagonal $\operatorname{U}(1)$. Define a corresponding meromorphic form on the torus $\mathbb{T}(d) = \left(\prod_{i\in Q_0}(\C^*)^{d_i}\right)/\C^*$, given explicitly by 
	\begin{align}\label{MeroForms}
		&\nonumber Z_{Q,d}(u) = \prod_{i\in Q_0}\prod^{d_i}_{s \neq s' = 1} \frac{u_{i,s'} - u_{i, s}}{u_{i, s} - u_{i,s'}-1} \\&\prod_{\{i\to j\} \in \bar{Q}_1}\prod^{d_i}_{s = 1}\prod^{d_j}_{s'=1}\left(\frac{u_{j,s'} - u_{i,s} + 1 - \frac{1}{2}R_{ij}}{u_{i,s} - u_{j,s'} + \frac{1}{2}R_{ij}}\right)^{\bra j, i\ket}  \bigwedge_{(i, s) \in Q_0 \times \{1,\ldots,d_i\} \setminus (\bar{i}, d_{\bar{i}})} du_{i,s},
	\end{align}
	for any fixed ordering of the pairs $(i, s) \in Q_0 \times \{1,\ldots,d_i\} \setminus (\bar{i}, d_{\bar{i}})$. Here $\bar{i} \in Q_0$ is a fixed reference vertex, {$\bar{Q}$ denotes the reduced quiver obtained by suppressing multiple arrows, and $\bra -, - \ket$ is the skew-symmetrised Euler form of $Q$.} The omission of a variable $u_{\bar{i}, N_{\bar{i}}}$ corresponds to quotienting the gauge group by the diagonal $\operatorname{U}(1)$. The variables $R_{ij}$ denote (a priori) arbitrary parameters, known as $R$-charges (analogous to $r_i$ in Example \ref{GLSMExm}). Fix a real stability vector $\zeta = \{\zeta_i,\,i \in Q_0\}$ (analogous to $t'$ in Example \ref{GLSMExm}). Then a physical conjecture (see \cite{beau}, Section 3), proved mathematically in many cases {(see \cite{ontani}, Theorem 5.1, \cite{OntaniStoppa}, Theorem 1.2)}, predicts the identity 
	\begin{equation}\label{DTJK}
	\bar{\chi}_Q(d, \zeta) = \frac{1}{d!}\jk_{\zeta}(Z_{Q, d}(u)) 
	\end{equation} 
{between a Jeffrey-Kirwan residue of a meromorphic form $Z_{Q, d}(u)$} (depending on the stability parameter $\zeta$) and a Donaldson-Thomas invariant $\bar{\chi}_Q(d, \zeta)$ of the quiver $Q$ {(as defined in great generality in \cite{MeinhardtReineke}, Definition 4.2)}. In particular, $\jk_{\zeta}(Z_Q(d), \zeta)$ should not depend on $R_{ij}$. 
	
\end{exm}

\subsection{Main result} Our aim here is to provide an analogue of Examples \ref{GLSMExm}, \ref{QuiverExm} in the context of general \emph{consistent scattering diagrams in} {$M_{\R} := M\otimes_{\Z}\R$}, where $M$ denotes a lattice of finite rank $r$, {as defined e.g. in \cite{leung}, Definitions 1.4 and 1.5 (we provide some background in Section \ref{MCSec}). Roughly speaking, such a scattering diagram $\fD$ is given by a collection of subsets (\emph{walls}) $P \subset M_{\R}$ decorated by automorphisms $\Theta_P$ of a certain complete graded $\C$-algebra attached to $P$ (\emph{wall-crossing automorphisms}), and the \emph{consistency} condition requires the vanishing of monodromy (composition of wall-crossing automorphisms) along loops. The most important application of consistent scattering diagrams (contained in general tropical affine manifolds with singularities) is the construction of mirror families, as in the work of Kontsevich and Soibelman \cite{KontsevichSoibelmanScattering} and in the Gross-Siebert programme (see e.g. \cite{GrossSiebertSympos}).     

As we explain in Section \ref{EnumerativeSubSec} below, our aim is motivated by the enumerative interpretations of such scattering diagrams. The idea that we develop in this paper is that formulae similar to those appearing in Examples \ref{GLSMExm}, \ref{QuiverExm},} involving Jeffrey-Kirwan residues, should exist quite generally for the coefficients of a consistent scattering diagram $\fD$, as an application of the remarkable approach developed by Chan-Leung-Ma and Leung-Ma-Young (see \cite{chan, leung}), using asymptotic analysis of Maurer-Cartan solutions in an appropriate differential graded Lie algebra (following proposals of Fukaya \cite{fukaya}). We refer to this approach as \emph{Maurer-Cartan scattering}.

So, in this paper, finite dimensional representations of a reductive group $G$ are replaced by an initial {(not necessarily consistent)} scattering diagram $\fD_{in} \subset M_{\R} \cong \R^r$, with values in a {Lie algebra $\fh$ satisfying certain assumptions spelled out in \cite{leung}, Section 1.1.1 (this is known as a \emph{tropical} Lie algebra). In particular, $\fh$ is graded} by $M^+_{\sigma} = (M \cap \sigma )\setminus \{0\}$ for a strictly convex cone $\sigma \subset M_{\R}$. 

{We consider the problem of computing the \emph{consistent completion} $\fD$ of $\fD_{in}$ in the sense of \cite{leung}, Theorem 1.6. This is a consistent scattering diagram obtained from $\fD_{in}$ by adding only \emph{outgoing} walls in the sense of \cite{leung}, Definition 1.3. The process of passing from $\fD_{in}$ to the consistent completion $\fD$ can be understood as eliminating all monodromy along loops and is uniquely defined up to a natural \emph{equivalence} relation on scattering diagrams, defined in \cite{leung}, Definition 1.5 (2) (see e.g. \cite{leung}, Theorem 1.6 for this result).

Thus we need to compute the coefficients of $\fD$, i.e. the elements $c_m \in \fh$ appearing as the graded pieces of the infinitesimal generator $\log\Theta_P$, where $\Theta_P$ denotes the automorphism attached to a wall $P$ of $\fD$.} 
 
In Section \ref{MCSec} we will recall the approach of \cite{chan, leung} in which the {consistent} completion $\fD$ is described in terms of (trivalent, labelled) trees $L$, which are essentially the possible {combinatorial} types of {certain embedded graphs, known as \emph{tropical discs},} contained in {the initial scattering diagram} $\fD_{in}$. {The boundaries (more precisely, parts of the boundaries known as the \emph{stops})} of such tropical discs with type $L$ trace a subset $P_{L} \subset M_{\R}$, a prospective wall of the consistent completion $\fD$. Also attached to $L$ is a degree $m_{L} \in M^+_{\sigma}$ and a ``multiplicity" $g_{L} \in \fh_{m_{L}}$.

This is all quite similar to other existing approaches (see e.g. \cite{fgs, gps, mandel}), but what is \emph{specific to Maurer-Cartan scattering} is the construction, for each type $L$, of a \emph{quadratic potential function} $W_{L}$ (defined on $M_{\R} \times \R^{|L^{[1]}|}$), with an associated Gaussian type bump $1$-form $\alpha_{L, \hbar}$, such that the contribution of $L$ to the consistent scattering diagram can be computed as the $\hbar \to 0$ limit of the \emph{transversal integral} $-\int_{\rho} \alpha_{L, \hbar}$ along various strata of $P$ ($\rho$ denoting a generic affine line intersecting transversely the interior of some stratum). In the present paper, these potential functions are what allows us to relate the consistent completion to residues of meromorphic forms. The quadratic property of $W_{L}$ is important for the relation to Jeffrey-Kirwan residues: we use that the first derivatives of the potential are affine linear and so define an arrangement of hyperplanes (in terms of the toric reductions studied in \cite{szenesvergne} they also define a system of weights).

In Section \ref{UnfoldingSubsec}, we adapt the above construction to \emph{unfoldings} $\widetilde{L}$ of $L$, which are types of tropical discs contained in an \emph{unfolded} diagram $\fD_{in, c}$, depending on a parameter $c \in (M_{\R})^J$ (for a suitable finite set $J$), see Definition \ref{UnfoldingDef}. As we explain, this unfolding is the \emph{minimal} deformation that allows one to achieve the required transversality. In particular, {our notion of unfolding is much simpler than the usual factorisation/deformation method for the computation of the consistent completion (described e.g. in \cite{gps}, Section 1.4 and \cite{leung}, Section 1.3), that is, it does not require to fix a particular factorisation of the wall-crossing automorphisms over an auxiliary ring.} 

{Let $\fD_{c}$ denote the consistent completion of the unfolding $\fD_{in, c}$,} see Definition \ref{UnfoldingComplDef}. By standard arguments, the consistent completions $\fD_{c}$ and $\fD$ are essentially equivalent as scattering diagrams (see Section \ref{UnfoldingSubsec} for a precise statement). 

In Section \ref{MainThmSubSec} we introduce a rational function $Z^{(P, m)}_{c, \delta}(s)$ of several complex variables $s = (s_1, \ldots, s_K)$, for sufficiently large $K$, with target a suitable space of  $\fh$-valued $1$-forms on $M_{\R}$. It is given by a finite sum 
\begin{align}\label{GlobalZIntro} 
Z^{(P, m)}_{c, \delta}(s) = \sum^K_{k = 1}\sum_{\widetilde{L} \in \cS_k(P, m)}\frac{2^{(k-1)/2}   g_{\widetilde{L}} (\pi \hbar)^{-1/2} e^{-W_{\widetilde{L}, c}(p)/\hbar} r_{\widetilde{L},c}}{|\Aut(\widetilde{L})||\del^2_s W_{L}|^{-1/2}\prod^{k-1}_{l = 1} \del_{s_l} W_{\widetilde{L},c+\delta_{\widetilde{L}}}(s)\prod^K_{j = k} s_j}
\end{align}
over all types $\widetilde{L}$ of tropical discs contained in $\fD_{in, c}$, with the same degree $m_{\widetilde{L}} = m \in M^+_{\sigma}$, and contributing nontrivially to the same wall $P$ of the consistent completion $\fD_{c}$ of $\fD_{in, c}$. Note that such $\widetilde{L}$ has an underlying type $L$ for $\fD_{in}$, and (since the potential is quadratic) the Hessian determinant $|\del^2_s W_{L}|$ is constant, both with respect to the complex variables $s$ and the parameters of $M_{\R}$. In the expression \eqref{GlobalZIntro}, $c$ is a sufficiently generic unfolding parameter, $p = \operatorname{Crit}(W_{\widetilde{L}, c})$ denotes the unique critical point of the potential $W_{\widetilde{L}, c}$, and $\delta = \{\delta_{\widetilde{L}}\}$ is a set of small and generic perturbation parameters. So, the factor $e^{-W_{\widetilde{L}, c}(p)/\hbar}$ is constant with respect to the complex variables $s$ and the product  $(\pi \hbar)^{-1/2} e^{-W_{\widetilde{L}, c}(p)/\hbar} r_{\widetilde{L},c}$ is a Gaussian type bump $1$-form on $M_{\R}$.  
\begin{thm}[Theorem \ref{MainThm}]\label{MainThmIntro} The coefficient $c_m$ of the consistent diagram $\fD_{c}$ along the wall $P$, in the fixed degree $m \in M^+_{\sigma}$, is given by the transversal integral of the Jeffrey-Kirwan residue $1$-form of $Z^{(P, m)}_{c, \delta}(s)$, namely
\begin{equation*}
c_m = -\lim_{\hbar \to 0}\int_{\rho} \jk^{\mathfrak{G}^{(P, m)}_{\delta}}_{\xi^{(P, m)}_{\delta}}(Z^{(P,m)}_{c, \delta}(s)),
\end{equation*}
where $\operatorname{JK}$ is taken with respect to the complex variables $s = (s_1, \ldots, s_K)$ and the set of affine hyperplanes
\begin{equation*}
\mathfrak{G}^{(P, m)}_{\delta}:= \{\del_{s_l} W_{\widetilde{L},c + \delta_{\widetilde{L}}}\!:\widetilde{L} \in \cS_k(P, m),\, l = 1,\ldots, k-1\} \bigcup \{s_j = 0,\,j = 1,\cdots, K\}, 
\end{equation*}
while the chamber $\xi^{(P, m)}_{\delta}$ is given by the positive linear span of $\mathfrak{G}^{(P, m)}_{\delta}$.
\end{thm}     
\begin{rmk} We will also prove a version of this result which does not depend explicitly on the $-\lim_{\hbar \to 0}\int_{\rho}$ operator, by passing to the \emph{theta functions} of the diagram $\fD_{c}$, see Theorem \ref{ThetaThm} below. 
\end{rmk}
\begin{rmk} There are at least some formal analogies between the rational function $Z^{(P, m)}_{c, \delta}(s)$ appearing in Theorem \ref{MainThmIntro} (constructed using Maurer-Cartan scattering) and the functions appearing in \eqref{QuasiMapZ} and \eqref{MeroForms}. Firstly, as already noted, the potentials $W_{\widetilde{L},c + \delta_{\widetilde{L}}}$ are quadratic and so the denominators involve affine linear functions $\del_{s_l} W_{\widetilde{L},c + \delta_{\widetilde{L}}}$, just as for \eqref{QuasiMapZ} and \eqref{MeroForms}. Moreover, in the special case when each graded piece $\fh_m \subset \fh$ is $1$-dimensional, with a distinguished basis vector $h_m$, we have an expression 
\begin{align*} 
& Z^{(P, m)}_{c, \delta}(s) = \\&\sum_{k > 0}\sum_{\widetilde{L} \in \cS_k(P, m)}\frac{2^{(k-1)/2} (\pi \hbar)^{-1/2} e^{-W_{\widetilde{L}, c}(p)/\hbar}r_{\widetilde{L},c}}{|\Aut(\widetilde{L})||\del^2_s W_{L}|^{-1/2}\prod^K_{j = k} s_j} \prod^{k-1}_{l = 1} \left(\frac{\del_{s_l} W_{\widetilde{L},c+\delta_{\widetilde{L}}}(s) + 1}{\del_{s_l} W_{\widetilde{L},c+\delta_{\widetilde{L}}}(s)}\right)^{\lambda_{\widetilde{L}, l}} h_{m_{L}}
\end{align*}
for unique scalars $\{\lambda_{\widetilde{L}, l}\}$, which is formally more similar to \eqref{MeroForms} (see Section \ref{MainThmSec}). This applies in particular to the important case of the Poisson algebra of formal Hamiltonian vector fields on an algebraic torus appearing e.g. in \cite{gps}, {Section 0.1}. 
\end{rmk}
\subsection{Relation to enumerative geometry}\label{EnumerativeSubSec} Consistent scattering diagrams in $M_{\R}$ are known to encode interesting enumerative invariants. There is a rather large literature on this topic, and we can only recall a few of these results here.
 
It was already noted in \cite{gp}, {Section 2.7, (based on Reineke's result \cite{Reineke}, Theorem 2.1),} that certain consistent scattering diagrams compute quiver invariants. {We will recall this is a special case in Example \ref{KroneckerExm}}. Much more generally, one can define a consistent scattering diagram using the invariants of any $2$-acyclic quiver with finite potential (see \cite{BridgelandScattering}, Theorem 1.1). 

{The Gross-Siebert programme for the construction of mirror families also leads to the enumerative interpretation of certain consistent scattering diagrams in $M_{\R}$ in terms of Gromov-Witten theory, see in particular \cite{ArguzGross}, Theorem 1.2, \cite{gps}, Theorem 0.1.}

In all these cases, Theorem \ref{MainThmIntro} provides, in principle, \emph{a representation of these enumerative (quiver and Gromov-Witten) invariants as the Jeffrey-Kirwan residues of a canonically defined family of rational functions}, thus giving an analogue of the representation formulae discussed in Examples \ref{GLSMExm}, \ref{QuiverExm}. {We will illustrate this in some simple quiver cases in Examples \ref{QuiverExm1}, \ref{QuiverExm2}.} However, we should stress that it is not at all clear how to compare directly the Jeffrey-Kirwan formula obtained in the case of quivers using our Theorem \ref{MainThmIntro} to the identity \eqref{DTJK}. 

\subsection{Theta functions} An advantage of the Maurer-Cartan scattering method of \cite{chan, leung} is its flexibility. For example, as shown in \cite{leung}, {Theorem 3.20}, it extends to give a description of the \emph{theta functions} of the consistent completion of a given initial scattering diagram {(see \cite{leung}, Section 1.2 for their definition). Such theta functions play a crucial role in the application of consistent scattering diagrams to mirror symmetry. }

In our case, we may apply these results to the consistent completion $\fD_c$ of a generic unfolding $\fD_{in, c}$. In Section \ref{ThetaProp}, given a pair of degrees $m, \mu \in M^+_{\sigma}$ and a generic point $Q \in M_{\R} \setminus  \fD_{c}$, we define rational functions $Z^{(Q, m, \mu)}_{c, \delta}(s)$, with values in $\fh$, analogous to \eqref{GlobalZIntro}. A variant of the proof of Theorem \ref{MainThmIntro} gives the following.
\begin{thm}[Proposition \ref{ThetaProp}]\label{ThetaThm} The degree $m$ theta function of $\fD_{c}$ can be expressed as  
\begin{equation*}
\widetilde{\vartheta}_m(Q) = \sum_{\mu \in M^+_{\sigma}} \jk^{\mathfrak{G}^{\mu}_{\delta}}_{\xi^{\mu}_{\delta}}(Z^{(Q, m, \mu)}_{c, \delta}(s)), 
\end{equation*}
for suitable hyperplane arrangements $\mathfrak{G}^{\mu}_{\delta}$ and chambers $\xi^{\mu}_{\delta}$.
\end{thm}
\begin{rmk} As pointed out to us by Riccardo Ontani, under suitable (restrictive) assumptions on $(Q, m, \mu)$, the hyperplane arrangement $\mathfrak{G}^{\mu}_{\delta}$ and chamber $\xi^{\mu}_{\delta}$ define the data of commuting actions of a torus $T$ and an auxiliary $\C^*$ on a linear space $V$, together with a regular linearisation $\xi$ with orbifold quotient $X := V /\!\!/_{\xi}T$. Then, by the general theory of \cite{szenesvergne}, one can construct an equivariant class $\alpha$ on $X$ such that
\begin{equation*}
\jk^{\mathfrak{G}^{\mu}_{\delta}}_{\xi^{\mu}_{\delta}}(Z^{(Q, m, \mu)}_{c, \delta}(s)) = \int_{X} \alpha := \int_{X^{\C^*}}  \frac{\alpha}{e^{\C^*}(N_{X^{\C^*}/X})}.     
\end{equation*} 
In light of this we can think of the expression appearing in Theorem \ref{ThetaThm} as formally similar to the ``moduli-theoretic" expression for theta functions given in \cite{BridgelandScattering}, Theorem 1.4 (see also \cite{cheung}, Chapter 7), in the special case when $\fD$ is the consistent scattering diagram constructed from a quiver (as in \cite{BridgelandScattering}, Theorem 1.1). 
\end{rmk} 
{
\subsection{Possible developments}
As we recalled, representations of enumerative invariants involving JK residues usually appear as localisation formulae, either infinite dimensional, in physics (see e.g. \cite{beau, beht, cordova}), or finite dimensional, in geometry (see e.g. \cite{KimUeda, ontani, szenesvergne} and their references). The aim of our results in the present paper, Theorems \ref{MainThmIntro}, \ref{ThetaThm}, is to show that similar JK residue formulae also appear in a different context, namely computing the consistent completion $\fD$ of an initial scattering diagram $\fD_{in}$. 

It would be interesting to understand if Theorem \ref{MainThmIntro} can also be explained, at least formally, as a suitable (infinite dimensional) localisation formula. 

Conversely, it seems interesting to ask if Theorem \ref{MainThmIntro} can be used to show that certain enumerative invariants, for example suitable infinite sets of quasimap invariants as in Example \ref{GLSMExm}, appear as the coefficients of a consistent scattering diagram (roughly, when their JK representations obtained from localisation can be deformed to the JK residues appearing in Theorem \ref{MainThmIntro}).

We should mention that a different result which relates JK residues directly to a special class of two-dimensional scattering diagrams is described in \cite{OntaniStoppa}, Theorem 1.4 $(ii)$. As explained in \cite{OntaniStoppa}, Section 1.3, these JK residues appear naturally in the computation of certain explicit flat sections of a meromorphic connection $\nabla_{\fD}$, which determines the consistent completion $\fD$. On the other hand, Theorem \ref{ThetaThm} in the present paper shows that similar JK residues appear when computing the theta functions of $\fD_{c}$, which by the general result \cite{leung}, Theorem 3.20 are smooth flat sections of a deformed differential $d_{\fD_{c}}$ acting on a complex of smooth forms (defined in a general context in \cite{leung}, Section 3.3), and also determine $\fD_c$ and so $\fD$ (as an application of \cite{leung}, Theorem 3.15). As envisaged in \cite{OntaniStoppa}, Remark 1.8, this gives a sort of ``de Rham" approach for the computation of (the theta functions of) $\fD$ based on the deformed differential $d_{\fD_{c}}$, as opposed to a ``Dolbeault" approach based on the meromorphic connection $\nabla_{\fD}$. Jeffrey-Kirwan residues appear naturally in both approaches. (Note that there is also a ``Betti" approach to the computation of theta functions based on broken lines, see \cite{leung}, Section 1.2. By \cite{leung}, Theorem 3.20 these actually define flat sections of the deformed differential).
}
\newline

\noindent{\textbf{Acknowledgements.}} We are grateful to Riccardo Ontani for some very helpful conversations related to the present work, and to the Referees for several important suggestions on improving the manuscript. This research was carried out in the framework of the project PRIN 2022BTA242 ``Geometry of algebraic structures: moduli, invariants, deformations".

\section{Maurer-Cartan scattering}\label{MCSec} 
Let $\fh$ denote a Lie algebra satisfying the assumptions of \cite{leung}, Section 1.1.1. We refer to $\fh$ as a \emph{tropical} Lie algebra. In particular, $\fh$ is graded by {$M^+_{\sigma} = (M \cap \sigma )\setminus \{0\}$ for a strictly convex cone $\sigma \subset M_{\R}$.} We use this grading to define completions. For elements $m \in M$, $n \in N:= M^{\vee}$ satisfying $\bra m, n\ket=0$, there is an abelian {graded} Lie subalgebra $\fh^{\parallel}_{m, n} := \bigoplus_{k > 0 }\fh_{k m,n} \subset \fh_{m}$. 

We follow the definitions of \emph{(consistent) scattering diagram in $M_{\R}$} set out in \cite{leung}, Definitions 1.4 and 1.5. \emph{Walls} $\opw \subset M_{\R}$ are described by data $(m, n, P, \Theta)$, {where} $m \in M$ is primitive, $\bra m, n \ket = 0$, $P \subset m_0 + n^{\perp}$ is a {convex rational polyhedral subset (the \emph{support} of the wall)} and $\Theta$ is an {element of $\exp(\hat{\fh}^{\parallel}_{m, n})$ (known as the \emph{wall-crossing automorphism} of $P$), $\hat{\fh}^{\parallel}_{m, n}$ denoting the completion of $\fh^{\parallel}_{m, n}$ with respect to the grading}.  

We fix an \emph{initial scattering diagram} $\fD_{in} = \{\opw_i = (m_i, n_i, P_i, \Theta_i)\}_{i \in I}$ and write
\begin{equation*}
\log \Theta_i = \sum_{j > 0} g_{j i}
\end{equation*}
for the decomposition into graded pieces. {Recall that, as discussed in the Introduction, there exists an essentially unique consistent completion $\fD$ of $\fD_{in}$, up to a natural \emph{equivalence} relation on scattering diagrams (spelled out in \cite{leung}, Definition 1.5 (2)). See for example \cite{leung}, Theorem 1.6 for this result.}
{  
\begin{exm}\label{KroneckerExm} Let us briefly recall one special example of the general relation between quiver representations and consistent scattering diagrams mentioned in the Introduction (see Section \ref{EnumerativeSubSec}). This is the case of the Kronecker quiver with $\kappa$ arrows $Q_{\kappa}$ and its framed representations, discussed in \cite{gp}, Section 2. Let us denote by $\M^{(1, 0), B}_{\kappa}(d a, d b)$, $\M^{(0, 1), F}_{\kappa}(d a, d b)$ the moduli spaces of stable framed representations of $Q_{\kappa}$ as in \cite{gp}, Section 2.4, and by 
\begin{align*}
& B_{(a, b)} = 1 + \sum_{d > 0} \chi(\M^{(1, 0), B}_{\kappa}(d a, d b)) x^{d a} y^{d b},\\ 
& F_{(a, b)} = 1 + \sum_{d > 0} \chi(\M^{(0, 1), F}_{\kappa}(d a, d b)) x^{d a} y^{d b},
\end{align*}  
the generating series for their topological Euler characteristics.

Choose $M := \Z^2$ endowed with the skew-symmetric bilinear form $\{ - , - \}$ induced by a multiple of the determinant, i.e. 
\begin{equation*}
\{ (m^1, m^2), ((m')^1, (m')^2)\} := \kappa\left(m^1 (m')^2 - m^2 (m')^1\right). 
\end{equation*}
Fix the strictly convex cone $\sigma = \R^2_{\geq 0} \subset \R^2 = M_{\R}$ and set $M_{\sigma} = M \cap \sigma$, $M^+_{\sigma} = (M \cap \sigma )\setminus \{0\}$. The monoid algebra $\fh := \C[M_{\sigma}]$ becomes a Lie algebra with bracket 
\begin{equation*}
\{x_m, x_m'\} := \{m, m'\} x_{m + m'} 
\end{equation*}
(this satisfies the assumptions of \cite{leung}, Section 1.1.1, i.e. it is tropical). Clearly, $\fh$ is graded by $M^+_{\sigma}$. Consider the initial scattering diagram
\begin{equation*}
\fD_{in} = \{\opw_i = (m_i, n_i, P_i, \Theta_i)\}_{i = 1, 2},
\end{equation*}
where
\begin{align*}
& m_1 = (1, 0),\, n_1 = \{ -, m_1\},\, P_1 = \R m_1, \Theta_1(x) = x,\, \Theta_1(y) = y (1 + x)^{\kappa}, \\
& m_2 = (0, 1),\, n_2 = \{ -, m_2\},\, P_2 = \R m_2, \Theta_2(x) = x(1+y)^{-\kappa},\,\Theta_2(y) = y,
\end{align*}
where we write $x = x_{m_1}$, $y = x_{m_2}$. Let $\fD$ denote the consistent completion of $\fD_{in}$. As noted in \cite{gp}, Section 2.7, a result of Reineke (see \cite{Reineke}, Theorem 2.1) shows that (the natural equivalence class of) $\fD$ can be represented as
\begin{equation*}
\fD = \{\opw_{(a, b)} = (m_{(a, b)}, n_{(a, b)}, P_{(a, b)}, \Theta_{(a, b)})\}_{(a, b) \in \Z^2_{\geq 0}}, 
\end{equation*}
where
\begin{equation*}
m_{(a, b)} = (a, b),\,n_{a, b} = \{- , m_{a, b}\},\,P_{(a, b)} = \R_{>0} m_{(a, b)},\, 
\end{equation*}
and we have
\begin{equation*}
\Theta_{(a, b)}(x) = x f^{-b}_{(a, b)},\, \Theta_{(a, b)}(y) = y f^a_{(a, b)},
\end{equation*}
where
\begin{equation*}
f_{(a, b)} = B^{\frac{\kappa}{a}}_{(a, b)} = F^{\frac{\kappa}{b}}_{(a, b)}. 
\end{equation*}
Thus, in this case, the consistent completion $\fD$ of $\fD_{in}$ computes the Euler characteristics of the moduli spaces of framed representations of $Q_{\kappa}$. 

Note that a direct computation (see e.g. \cite{FilippiniStoppa}, Lemma 3.4) shows that the infinitesimal generators of the wall-crossing automorphisms in the initial scattering diagram are given by
\begin{equation*}
\log \Theta_i = \operatorname{Li}_2(-x_{m_i}) = \sum_{j > 0} \frac{(-1)^j x_{j m_i}}{j^2},\, i = 1, 2.
\end{equation*}
So we have
\begin{equation*}
g_{ji} = \frac{(-1)^j x_{j m_i}}{j^2},\, i = 1, 2. 
\end{equation*}
We will use this in Examples \ref{QuiverExm1}, \ref{QuiverExm2}.
\end{exm}
}
The main results of \cite{leung} are expressed in terms of sums over various sets of trees. These are all versions of \emph{directed, trivalent} trees $T$, with a decomposition of the set of vertices 
\begin{equation*}
\bar{T}^{[0]} = T^{[0]}_{in} \cup T^{[0]} \cup \{v_{out}\} 
\end{equation*}
(see \cite{leung}, Definition 1.11). Such a tree with $|T^{[0]}_{in}| = k$ is called a \emph{$k$-tree}.  

\emph{Labelled $k$-trees} $L \in \LT_k$ have a labelling of each incoming edge $e \in L^{[1]}_{in}$ by a wall $\opw_{i_e} = (m_{i_e}, n_{i_e}, P_{i_e}, \Theta_{i_e})$ of $\fD_{in}$ together with an element $m_{e} \in M^+_{\sigma}$ such that $m_e$ is a  multiple of the primitive element $m_{i_e}$. This induces a labelling of all edges of $L$ by elements of $M^+_{\sigma}$, by the rule $m_{e_3} = m_{e_1} + m_{e_2}$ , where $e_1, e_2$ are incoming at a vertex $v$ and $e_3$ is outgoing from $v$ {(resembling the balancing condition for tropical curves, see \cite{gps}, Definition 2.1 (2)).} The labelling $m_{e_{out}}$ is also denoted by $m_{L}$.   

As explained in \cite{leung}, Section 1.3 {(see especially the discussion immediately following Definition 1.16)}, $L \in \LT_k$ defines a \emph{tropical subspace} $P_{L} \subset M_{\R}$, given by the locus traced by the \emph{stops} {(i.e. the outgoing part of the boundary)} of the tropical discs of type $L$ contained in $M_{\R}$ {(\cite{leung}, Definition 1.16)}. 

\emph{Labelled ribbon $k$-trees} $\cL \in \LR_k$ carry the further datum of a cyclic order of the edges incident at $v$ for all vertices $v$ (this includes a cyclic order of the leaves $\cL^{[0]}_{in}$). Of course, each $\cL \in \LR_k$ has an underlying $L \in \LT_k$.

Using the labelled ribbon structure, one can attach naturally to $\cL \in \LR_k$ a pair $(n_{\cL}, g_{\cL}) \in N \times \fh^{\parallel}_{m_{\cL}, n_{\cL}}$ as in \cite{leung}, Definition 1.14, i.e. $n_{\cL}$ is determined by the tropical rule, $g_{\cL}$ by iterated Lie brackets, and the sign is fixed by the ribbon structure. In fact such a pair $(n_e, g_e)$ is canonically attached to each $e \in \cL^{[1]}$, and we have $(n_{\cL}, g_{\cL}) = (n_{e_{out}}, g_{e_{out}})$.    

The main object of study in \cite{leung} are certain special solutions $\Phi$ of the Maurer-Cartan (MC) equation
\begin{equation*}
d\Phi + \frac{1}{2}[\Phi, \Phi] = 0
\end{equation*}
with values in the completion $\hat{\calH}^*$ of the differential graded Lie algebra
\begin{equation*}
\calH^* = \calH^*(M_{\R}) = \bigoplus_{m \in M^+_{\sigma}} \big(\cW^0_*(M_{\R}) / \cW^{-1}_*(M_{\R})\big) \otimes_{\C} \fh_m.  
\end{equation*}  
Here, $\cW^0_*(M_{\R}) / \cW^{-1}_*(M_{\R})$ is a suitable functional space of differential forms on $M_{\R}$ (of all degrees), whose support concentrates along a tropical subspace $P \subset M_{\R}$ as $\hbar \to 0$, modulo forms that vanish as $\hbar \to 0$ (see \cite{leung}, Section 2.2.1). { The differential $d$ acts on the form part, while the Lie bracket on $\calH^*$ is a combination of the wedge product and the Lie bracket on $\fh$. These are well defined on $\hat{\calH}^*$ by \cite{leung}, Lemma 2.3.}  

These special solutions are constructed in { \cite{leung}, Section 2.3, by using a method originally due to Kuranishi \cite{Kuranishi} and explained in detail in \cite{chan}, Section 5.1. We note that, in general, Kuranishi's method involves the choice of a \emph{homotopy operator} $\opH$, acting on differential forms, with degree $-1$. In our present context, the choice of $\opH$ is explained in \cite{leung}, Section 3.1.1, and we will recall it briefly below.} 

{Thus, following this general method, solutions are expressed as a sum of $1$-forms indexed by labelled ribbon trees, of the form
\begin{equation}\label{KuranishiMCelement}
\Phi = \sum_{k > 0} \frac{1}{2^{k-1}} \sum_{\cL \in \LR_k} \opL_{k, \cL}(\Pi, \ldots, \Pi).
\end{equation}

Here the \emph{initial datum} $\Pi$ is given in terms of defining affine linear functions $\eta_i$ for $P_i$,  
\begin{equation*}
\Pi = \sum_{i \in I} -\delta_{P_i} \log(\Theta_i),\,\delta_{P_i} := \left(\frac{1}{\pi \hbar}\right)^{1/2} e^{-\eta^2_i/\hbar} d\eta_i.
\end{equation*}	

We note that, by construction, any internal or final edge of $\cL \in \LR_k$, thought of as an outgoing edge, defines an element of $M_{\R}$. The \emph{homotopy operator} $\opH$, placed at such an edge, acts on differential forms, with degree $-1$, by integrating along the linear flow in the direction of $M_{\R}$ corresponding to the outgoing edge (see \cite{leung} Section 3.1.1). 

Finally, the operation $\opL_{k, \cL}$, acting on $k$-tuples of elements of $\hat{\calH}^*$ and yielding a $1$-form, is defined in \cite{leung}, Definition 2.12. In brief, it is given by compositions of the Lie bracket, performed at incoming edges of $\cL$, and the (negative) homotopy operator $-\opH$, performed at outgoing edges, starting from the $M^+_{\sigma}$-graded components of the inputs, with respect to the gradings specified by the labelling of the leaves of $\cL$. 

We illustrate this definition concretely in Examples \ref{ResidueExm1} and \ref{ResidueExm2} below.}

By construction, we may write
\begin{equation*}
\opL_{k, \cL}(\Pi, \ldots, \Pi) = \alpha_{\cL} g_{\cL},
\end{equation*} 
where the $1$-form $\alpha_{\cL}$ is defined {by the Kuranishi method as above}, and \emph{can be written explicitly} (see \eqref{1formAlpha} in Section \ref{MCSec} below).  

\begin{rmk}\label{SignAmbiRmk} One shows that the value of $\opL_{k, \cL}(\Pi, \ldots, \Pi)$ does not depend on the ribbon structure of $\cL$, so both $\opL_{k, L}(\Pi, \ldots, \Pi)$ and the product $\alpha_{L} g_{L}$ are well defined for $L \in \LT_k$. On the other hand, we can also define each of $\alpha_{L}$ and $g_{L}$ in the same way, up to the \emph{same} choice of sign.  
\end{rmk} 

According to \cite{leung}, Lemma 3.10 and Lemma 3.11, there exists a distinguished polyhedral decomposition $\cP_{L}$ of the tropical subspace $P_{L}$ into cells $\sigma$, such that the transversal integral $\int_{\rho}\alpha_{L}$ is constant along the relative interior of each top dimensional cell $\sigma \in \cP^{[r-1]}_{L}$ (where $\rho$ denotes a generic affine line intersecting the relative interior of $\sigma$ transversely). 

Following \cite{leung}, Definition 3.13, one defines a new scattering diagram as the collection of walls
\begin{equation*}
\fD = \{(m_{\cL}, n_{\cL}, \sigma \subset P_{L}, \exp\big(-\frac{1}{|\Aut(L)|} \int_{\rho}\alpha_{L} g_{L}\big))\}.   
\end{equation*}  
{Note that, by construction, $\fD$ contains the initial scattering diagram $\fD_{in}$.

We can now state one of the main results of Leung, Ma and Young in \cite{leung}, building on the previous work of Chan, Leung and Ma \cite{chan}. 
\begin{thm}[\cite{leung}, Theorem 3.12 and Proposition 3.14]\label{LMYThm} $\fD$ is consistent.
\end{thm}
}
Naturally, the Maurer-Cartan equation plays a crucial role in the proof of consistency.
\subsection{Unfolding}\label{UnfoldingSubsec}
Let $\fD_{in} = \{\opw_i = (m_i, n_i, P_i, \Theta_i)\}_{i \in I}$ be an initial scattering diagram. We introduce a tautological notion of \emph{unfolding} of $\fD_{in}$, which we use to achieve certain transversality conditions in our applications of Maurer-Cartan scattering.

We work over the ring  
\begin{equation*}
R_{J} := \C[t_{i, j},\,i\in I,\, j \in J \subset \Z_{>0}]/(t^2_{i,j}), 
\end{equation*} 
where $J$ is a finite interval $\{1,\ldots, |J|\}$, as large as required. 
{
\begin{definition}\label{UnfoldingDef} We define the unfolding of $\fD_{in}$ with parameter $c \in (M_{\R})^{J}$ as the scattering diagram $\fD_{in, c}$ over $R_J$ with walls
\begin{equation*}
\opw_{i, j} = (m_i, n_i, P_{i, j} := P_i + c_{i,j}, \log \Theta_{i, j} := t_{i, j} \log \Theta_i).
\end{equation*}
\end{definition}\label{UnfoldingComplDef}
\begin{definition} We denote by $\fD_{c}$ the consistent completion of the unfolding $\fD_{in, c}$ (in the sense of \cite{leung}, Theorem 1.6).
\end{definition} 
}
Fix a wall $\opw \subset M_{\R}$ of $\fD_{c}$. The corresponding weight function $\log \Theta $ can be written uniquely as a polynomial in the variables $t_{i, j}$. We may define a new scattering diagram $\fD'_{c}$ by replacing each $\log \Theta$ with the evaluation of the corresponding polynomial at 
\begin{equation*}
t_{i, j} = \frac{1}{|J|}. 
\end{equation*}
Then, by construction, the asymptotic diagram $(\fD'_{c})_{as}$ (defined as in \cite{gps}, Definition 1.7, using essentially that $\fD'_{c}$ is contained in $M_{\R}$) is equivalent to the consistent completion of $\fD_{in}$, up to any fixed, arbitrary order, determined by $|J|$.
\section{MC scattering and (higher) residue pairing}\label{PairingSec} 

According to \cite{leung}, Lemma 3.4 and Lemma 3.9, we may write
\begin{equation*}
\opL_{k, \cL}(\Pi, \ldots, \Pi) = \alpha_{\cL} g_{\cL},
\end{equation*}
where $g_{\cL}$ is the combinatorial factor attached to $\cL$, depending only on $\cL$ and $\fh$ (but not on the Maurer-Cartan equation), while the $1$-form part is given by 
\begin{equation}\label{1formAlpha}
\alpha_{\cL}(x) = \left(\frac{1}{\pi \hbar}\right)^{k/2}\int_{[-\infty, 0]^{\cL^{[1]}}\times \{x\}} e^{-\big(\sum^k_{j=1} (\tau^{\fe_j})^* \eta^2_{i_{e_j}}\big)/\hbar} \vec{\tau}^*\big( d\eta_{e_1} \wedge \cdots \wedge d\eta_{e_k} \big).
\end{equation}
Here, the pullback operations $(\tau^{\fe_j})^*$, $\vec{\tau}^*$, defined in \cite{leung}, Section 3.12, are simply those induced by the linear flows which appear when applying the (negative) homotopy operator $-\opH$ in the Kuranishi method (see Examples \ref{ResidueExm1} and \ref{ResidueExm2} below for more details).

\emph{In the following, we will regard the right hand side of \eqref{1formAlpha} as a complex oscillatory integral, so that its $\hbar \to 0$ asymptotics can be expressed in terms of the (higher) residue pairing, e.g. by the results of \cite{CoatesCortiIritani_hodge}, Section 6.} 

We will always assume the nondegeneracy condition 
\begin{equation*}
n_{\cL} \neq 0.
\end{equation*}
This is not restrictive since by \cite{leung}, Lemma 3.4 (1), if $n_{\cL} = 0$ then $\alpha_{\cL} = 0$.

We introduce the complex algebraic torus 
\begin{equation*}
\TT_{\cL}= (\C^*)^{\cL^{[1]}} \cong (\C^*)^{k-1},  
\end{equation*}
and denote the natural complex coordinates and holomorphic volume form on its Lie algebra $\Lie(\TT_{\cL})$ by $s = (s_1, \ldots, s_{k-1})$ and  
\begin{equation*}
\Omega_{\cL} = \bigwedge^{k-1}_{j = 1}  d s_j  
\end{equation*}
(we use here that the ribbon structure on $\cL$ induces a cyclic order of $\cL^{[1]}$). Note that $(\tau^{\fe_j})^* \eta^2_{i_{e_j}},\,j =1,\ldots, k$ and $\vec{\tau}^*\big( d\eta_{e_1} \wedge \cdots \wedge d\eta_{e_k} \big)$ have \emph{canonical  holomorphic extensions} to $\TT_{\cL}$, since the functions $\eta_{i_{e_j}}$ are affine linear (over $\R$) and $(\tau^{\fe_j})^*$, $\vec{\tau}^*$ are induced by linear flows (over $\R$).

We also upgrade $\hbar$ and $x$ to complex variables $-z$ and $u = (u_1, \ldots, u_r)$. Then, setting 
\begin{align*}
& W_{\cL}(s) =  \sum^k_{j=1} (\tau^{\fe_j})^* \eta^2_{i_{e_j}},\, \Gamma = [-\infty, 0]^{\cL^{[1]}} \subset \Lie(\TT_{\cL}),\\
& r_{\cL}(s)= \frac{\vec{\tau}^*\big( d\eta_{e_1} \wedge \cdots \wedge d\eta_{e_k} \big)}{\Omega_{\cL}},     
\end{align*}
we extend \eqref{1formAlpha} to the complex oscillatory integral
\begin{equation}\label{cplx1formAlpha}
\alpha_{\cL}(u) = \left(\frac{1}{\pi(-z)}\right)^{k/2}\int_{\Gamma} e^{W_{\cL}(s)/z} r_{\cL}(s)\Omega_{\cL}\big|_u.
\end{equation}

\begin{rmk}
By construction, $r_{\cL}$ is a $1$-form on $M_{\R}$. It is important to keep in mind that here (and in all the following), all our quantities are valued in smooth $1$-forms in the variables $u_1, \ldots, u_{r-1}$.
\end{rmk}
\begin{rmk} Following Remark \ref{SignAmbiRmk}, one has a similar integral yielding the $1$-form $\alpha_{L}(u)$, with the only difference that the factor $r_{L}$ is only well-defined up to sign. The product $r_{L} g_{L}$, on the other hand, is well-defined.
\end{rmk}
\subsection{Behaviour of the oscillatory integral} The potential $W_{\cL}(s)$ is a quadratic function. By \cite{leung}, Lemma 3.3 the nondegeneracy condition $n_{\cL} \neq 0$ implies that, \emph{for fixed $u = (u_1, \ldots, u_r) \in M_{\R}$}, there is a unique critical point $p$ given by the solution of the linear system
\begin{equation*}
\nabla_s W_{\cL}|_p = 2\sum^k_{j=1} (\tau^{\fe_j})^* \eta_{i_{e_j}} \frac{\del}{\del s_l}\big(\tau^{\fe_j}\big)|_p = 0,\,l=1,\ldots,k-1.
\end{equation*}
Moreover, by construction, taking the unique critical point $p$, \emph{as $u \in M_{\R}$ varies}, induces a linear map
\begin{equation*}
p_{\cL}\!: M_{\R} \to \Lie(\TT_{\cL}).  
\end{equation*}
The asymptotic behaviour of the integral \eqref{cplx1formAlpha} depends on whether
\begin{enumerate}
\item[(i)] $p \in \Int_{re}(\Gamma) = \Gamma\setminus\del\Gamma$,\,
\item[(ii)] $p \in \del\Gamma$, or
\item[(iii)] $p \in \Lie(\TT_{\cL}) \setminus \Gamma$.
\end{enumerate}
\subsubsection{The case $p \in \Int_{re}(\Gamma) = \Gamma\setminus\del\Gamma$} Recall that, according to \cite{CoatesCortiIritani_hodge}, Definition 6.6, the \emph{higher residue pairing} (with respect to the potential $W_{\cL}$) is given, for arbitrary $f_1, f_2$, by 
\begin{equation*}
P(f_1 \Omega_{\cL}, f_2 \Omega_{\cL}) = \overline{\Asymp_p(e^{W_{\cL}/z} f_1 \Omega_{\cL})} \Asymp_p(e^{W_{\cL}/z} f_2\Omega_{\cL}),
\end{equation*}
where the \emph{asymptotic expansion operator} $\Asymp_p$ is defined by requiring, for all $f$, 
\begin{equation*}
\int_{\Gamma} e^{W_{\cL}/z} f \Omega_{\cL} \sim e^{W_{\cL}(p)/z} (-2\pi z)^{(k-1)/2} \Asymp_p(e^{W_{\cL}/z} f \Omega_{\cL}),
\end{equation*}  
and we set
\begin{equation*}
\overline{\Asymp_p(e^{W_{\cL}/z} f \Omega_{\cL})} = \Asymp_p(e^{W_{\cL}/z} f \Omega_{\cL})|_{z \to -z}. 
\end{equation*}
Note that we have 
\begin{equation*}
\Asymp_p(e^{W_{\cL}/z} f \Omega_{\cL}) = \frac{1}{\sqrt{\det \nabla^2 W_{\cL}(p)}}\big(f(p) + a_1z + a_2 z^2 + \cdots\big). 
\end{equation*} 
In particular, the first coefficient in the $z$ expansion of $P(f_1 \Omega_{\cL}, f_2 \Omega_{\cL})$ coincides with the {\emph{Grothendieck residue pairing}} as defined e.g. in \cite{GriffithsHarris}, p. 659 (see also \cite{Hertling}, Section 10.4), namely
\begin{align*}
\res_{W_{\cL}}(f_1, f_2) &= \Res_{\{p\}}\left(\left(\frac{1}{2\pi \ii}\right)^{k-1} \frac{f_1 f_2 \Omega_{\cL}}{\prod^{k-1}_{l = 1} \del_{s_l} W_{\cL}} \right)\\
&= \left(\frac{1}{2\pi \ii}\right)^{k-1} \int_{C_{\varepsilon}} \frac{f_1 f_2 \Omega_{\cL}}{\prod^{k-1}_{l = 1} \del_{s_l} W_{\cL}} ,
\end{align*} 
where the integration cycle is given by
\begin{equation*}
C_{\varepsilon} = \{ |\del_{s_l} W_{\cL}|  = \varepsilon\},
\end{equation*}
with positive orientation with respect to the ordered set $(\del_{s_l} W_{\cL}), \,l = 1,\ldots, k-1$. 

Applying this to our case, we find the leading order $z \to 0$ term
\begin{align*}
&\alpha_{\cL}(u) \sim\\
& \left(\frac{1}{\pi(-z)}\right)^{k/2} \frac{e^{W_{\cL}(p)/z}}{|\del^2_s W_{\cL}|^{-1/2}} (-2\pi z)^{(k-1)/2} \left(\frac{1}{2\pi \ii}\right)^{k-1} \int_{|\del_{s_l} W_{\cL}|  = \varepsilon} \frac{r_{\cL} \Omega_{\cL}}{\prod^{k-1}_{l = 1} \del_{s_l} W_{\cL}},
\end{align*}
where we write $|\del^2_s W_{\cL}|$ for the \emph{constant} Hessian determinant. We have thus established the following result.
\begin{prop}\label{LeadingOrderProp} The leading order term in the $z \to 0$ asymptotic expansion of the $1$-form $\alpha_{\cL}$ along the locus $p^{-1}_{\cL}(\Int_{re}(\Gamma))$ is given by 
\begin{align*}
 2^{(k-1)/2} (\pi(-z))^{-1/2} \frac{e^{W_{\cL}(p)/z}}{|\del^2_s W_{\cL}|^{-1/2}} \left(\frac{1}{2\pi \ii}\right)^{k-1} \int_{C_{\varepsilon}} \frac{r_{\cL} \Omega_{\cL}}{\prod^{k-1}_{l = 1} \del_{s_l} W_{\cL}}  
\end{align*}
for all sufficiently small $\varepsilon > 0$, where each term $\del_{s_l} W_{\cL}$ is given by an \emph{affine linear} function 
\begin{equation*}
\del_{s_l}W_{\cL} = 2\sum^k_{j=1} (\tau^{\fe_j})^* \eta_{i_{e_j}} \frac{\del}{\del s_l}\big(\tau^{\fe_j}\big).
\end{equation*}
\end{prop}
Indeed, by construction, each term $\del_{s_l}W_{\cL}$ is an affine linear function in the variables $s_1, \cdots, s_{k-1}$, with constant term given by a linear function in the variables $u_1, \ldots, u_r$. As above, assuming the nondegeneracy condition $n_{\cL} \neq 0$, \cite{leung}, Lemma 3.3 implies that the linear span of the affine functions $\del_{s_l}W_{\cL},\,l=1,\ldots, k-1$ has full rank. Setting 
\begin{equation*}
\zeta_l = \del_{s_l}W_{\cL},\,l=1,\ldots, k-1, 
\end{equation*} 
we write 
\begin{equation*}
\frac{r_{\cL} \Omega_{\cL}}{\prod^{k-1}_{l = 1} \del_{s_l} W_{\cL}} = Z_{\cL}(\zeta) d\zeta_1 \wedge \cdots \wedge d\zeta_{k-1},
\end{equation*}
where $Z_{\cL}(\zeta)$ is a meromorphic function, defined in a neighbourhood of $0 \in (\C^*)^{k-1}$, \emph{valued in smooth $1$-forms in the variables $u_1,\ldots,u_r$}. 

Then, {by the standard residue theorem}, we can express the result of Proposition \ref{LeadingOrderProp}, {which expresses the leading asymptotics as an iterated integral}, alternatively in terms of an \emph{iterated residue} of a meromorphic function, {in the sense recalled e.g. in \cite{OntaniStoppa}, Appendix A.}
\begin{cor}\label{LeadingOrderCor} The leading order term in the $z \to 0$ asymptotic expansion of the $1$-form $\alpha_{\cL}$ along the locus $p^{-1}_{\cL}(\Int_{re}(\Gamma))$ is given by 
\begin{equation*}
 2^{(k-1)/2} (\pi(-z))^{-1/2} \frac{e^{W_{\cL}(p)/z}}{|\del^2_s W_{\cL}|^{-1/2}} \IR_0(Z_{\cL}(\zeta)),
\end{equation*}
{where $\IR_0(Z_{\cL}(\zeta))$ denotes the iterated residue of the meromorphic function $Z_{\cL}(\zeta)$, computed with respect to the ordered variables $\zeta_1, \ldots, \zeta_{k-1}$ at the point $0 \in (\C^*)^{k-1}$ (see e.g. \cite{OntaniStoppa}, Appendix A).}
\end{cor}
\begin{exm}\label{ResidueExm1} Recall the $1$-form $\alpha_{\cL} g_{\cL}$ is defined by Kuranishi's method, applying iteratively the Lie bracket in $\calH^*$ followed by the (negative) homotopy operator $-\opH$ along the labelled ribbon tree $\cL$. In particular restricting to coefficients in $\cW^0_*(M_{\R}) / \cW^{-1}_*(M_{\R})$ gives the $1$-form $\alpha_{\cL}$. 

Fix $M = \Z^2$ and consider $\cL \in \LR_2$ {with roots decorated by $\{m_1 = (1, 0), m_2 = (0,1)\}$}. Then we have  
\begin{equation*}
\alpha_{\cL} = -\opH(\delta \wedge \eta),
\end{equation*} 
where
\begin{equation*}
\delta = \left(\frac{1}{\pi \hbar}\right)^{1/2} e^{-y^2/\hbar} dy,\,\eta = \left(\frac{1}{\pi \hbar}\right)^{1/2} e^{-x^2/\hbar} dx
\end{equation*}
and the linear flow corresponding to $\opH$ has infinitesimal generator
\begin{equation*}
\frac{\del}{\del s} = (-1, -1) \in M.
\end{equation*}
We have by definition
\begin{equation*}
\opH(\delta \wedge \eta) = \left(\frac{1}{\pi \hbar}\int^0_{-\infty} e^{-((x-s)^2+(y-s)^2)/\hbar} ds \right)\iota_{\frac{\del}{\del s}} (dy\wedge dx), 
\end{equation*}
and so
\begin{equation*}
\alpha_{\cL}(u) = \left(\frac{1}{\pi (-z)}\int^0_{-\infty} e^{((u_1-s)^2+(u_2-s)^2)/z} ds \right) d(u_1-u_2). 
\end{equation*}
Thus, 
\begin{align*}
& \cW_{\cL} = (u_1-s)^2+(u_2-s)^2,\,\del_s \cW_{\cL} = -2(u_1+u_2 - 2s),\\
& p = \frac{u_1+u_2}{2},\,\cW_{\cL}(p) =  \frac{(u_1-u_2)^2}{2} 
\end{align*}
and the leading order term in the $z \to 0$ asymptotic expansion of $\alpha_{\cL}(u)$ along the locus
\begin{equation*}
p^{-1}_{\cL}(\Int_{re}(\Gamma)) = \{u_1 + u_2 < 0\} \subset M_{\R}
\end{equation*}
is given by
\begin{equation*}
\frac{1}{2\pi \ii}  \int_{C_{\varepsilon}} \frac{1}{-2(u_1+u_2 - 2s)} ds \cdot 2^{1/2}(\pi(-z))^{-1/2} 2 e^{\frac{(u_1-u_2)^2}{2 z}} d(u_1-u_2),
\end{equation*}
where 
\begin{equation*}
C_{\varepsilon} = \{s : |u_1 + u_2 - 2s| = \varepsilon/2\} 
\end{equation*}
with positive orientation. Setting $\zeta = -2(u_1+u_2 - 2s)$, we have 
\begin{equation*}
Z_{\cL}(\zeta) = \left( \frac{1}{4 \zeta}\right) d(u_1-u_2),\,\alpha_{\cL} \sim 2^{1/2}(\pi(-z))^{-1/2} 2 e^{\frac{(u_1-u_2)^2}{2 z}} \IR_0(Z_{\cL}(\zeta)).
\end{equation*}
\end{exm}
\subsubsection{The case $p \in \del\Gamma$} By definition, the boundary $\del\Gamma$ is the real part of the locus
\begin{equation*}
\{\prod^{k-1}_{l=1} s_l = 0\} \subset \Lie(\TT_{\cL}).
\end{equation*} 
Thus, $\del\Gamma$ admits a natural stratification, and the asymptotic expansion of the oscillatory integral \eqref{cplx1formAlpha} for $p$ in each fixed stratum is of the form 
\begin{equation*}
 2^{(k-1)/2} \kappa_{\cL} (\pi(-z))^{-1/2} e^{W_{\cL}(p)/z} 
\end{equation*}
for a suitable constant $\kappa_{\cL}$.
\begin{exm}\label{ResidueExm2} We consider $M = \Z^2$ and $\cL \in \LR_3$ {with roots decorated by $\{m_1 = (1,0), m_2 = (0,1), m_3 = (0,1)\}$. Then} 
\begin{equation*}
\alpha_{\cL} = -\opH(-\opH(\delta \wedge \eta) \wedge \eta) 
\end{equation*}
{where} the linear flow corresponding to the first insertion of $\opH$ has infinitesimal generator
\begin{equation*}
\frac{\del}{\del s_2} = (-1, -2).
\end{equation*}
By our previous computations, the right hand side is given by
\begin{align*}
 & \left(\frac{1}{\pi \hbar}\right)^{3/2} \left(\int_{(-\infty,0]^2} e^{-((x-s_2)^2+(x-s_1-s_2)^2+(y-s_1-2s_2)^2)/\hbar} ds_1 ds_2 \right)\\&\otimes \iota_{\frac{\del}{\del s_2}}((dx - dy )\wedge dx),
\end{align*}
so we have
\begin{equation*}
\alpha_{\cL}(u) = \left(\frac{1}{\pi (-z)}\right)^{3/2} \int_{(-\infty,0]^2} e^{W_{\cL}/z} ds_1 ds_2\,d(2u_1-u_2),
\end{equation*}
where
\begin{equation*}
W_{\cL} = (u_1-s_2)^2+(u_1-s_1-s_2)^2+(u_2-s_1-2s_2)^2.
\end{equation*}
We compute  
\begin{align*}
& p = (0, \frac{u_1+u_2}{3}),\, \cW_{\cL}(p) = \frac{1}{3}(2u_1-u_2)^2,\\
& \del_{s_1} \cW_{\cL} = 2 (2 s_1+3 s_2-u_1-u_2),\,\del_{s_2} \cW_{\cL} = 2 (3 s_1 + 6  s_2 -2 u_1 -2 u_2). 
\end{align*}
So, in this case, the critical point $p(u)$ \emph{always lies in the stratum $s_1 = 0$ of $\del \Gamma$}, and it lies in the interior of this stratum iff $u_1 + u_2 < 0$. 
\end{exm}
\subsubsection{The case $p \in \Lie(\TT_{\cL}) \setminus \Gamma$} In this case the transversal integrals $\int_{\rho}\alpha_{\cL}$ vanish in the $\hbar \to 0$ limit, and there is no contribution to the consistent scattering diagram.
\section{Unfolding}\label{UnfoldingSec}
Here we examine the effect of unfolding the initial scattering diagram $\fD_{in}$ (as described in Section \ref{UnfoldingSubsec}) on the oscillatory integrals \eqref{cplx1formAlpha}. 

By construction, this means that we apply Maurer-Cartan scattering to the unfolded scattering diagram $\fD_{in, c}$ over the ring $R_J$. Let $\widetilde{\cL}$ be a labelled ribbon tree for the unfolded diagram $\fD_{in, c}$ (with underlying labelled tree $\widetilde{L}$ for $\fD_{in, c}$). Then, there is a unique labelled ribbon tree $\cL$ \emph{for the original diagram $\fD_{in}$} such that 
\begin{equation}\label{IntegralAlphaUnfolded} 
\alpha_{\widetilde{\cL}}(x) = \left(\frac{1}{\pi \hbar}\right)^{k/2}\int_{\Gamma_{\cL}} e^{-\big(\sum^k_{j=1} (\tau^{\fe_j})^* \eta^2_{i_{e_j}, c}\big)/\hbar} \vec{\tau}^*\big( d\eta_{e_1,c} \wedge \cdots \wedge d\eta_{e_k, c} \big), 
\end{equation}
where $\Gamma_{\cL} = [-\infty, 0]^{\cL^{[1]}}\times \{x\}$ and the new defining affine functions are determined by 
\begin{equation*}
\eta_{P_{i, j}, c}(x):= \eta_{P_i}(x - c_{i,j}).
\end{equation*}

Note the following crucial property of the unfolding procedure: \emph{for each fixed labelled tree $\widetilde{\cL} \in \LR_k(\fD_{in, c})$, by the relation $t^2_{i, j} = 0$ in $R_J$, the shifts $c_{i, j}$ can be chosen independently for $j = 1, \ldots, k$}. 

Recall that, as $u \in M_{\R}$ varies, taking the unique critical point $p(u)$ of $W_{\cL}$ induces a linear map $p_{\cL}\!: M_{\R} \to \Lie(\TT_{\cL})$. In the unfolded diagram, for each fixed labelled ribbon tree $\widetilde{\cL}$ for $\fD_{in, c}$, this is enhanced to a map 
\begin{equation*}
p_{\widetilde{\cL}}\!: M_{\R} \times (M_{\R})^{J} \to \Lie(\TT_{\widetilde{\cL}}), 
\end{equation*}
where $p_{\widetilde{\cL}}$ is the unique critical point of the quadratic potential
\begin{equation*}
W_{\widetilde{\cL}, c} = \sum^k_{j=1} (\tau^{\fe_j})^* \eta^2_{i_{e_j}, c}
\end{equation*}
with respect to the variables $s_1, \ldots, s_{k-1}$, as a function of $x$ and $c$. 

By the property spelled out above, for all fixed $x$, the map $p_{\widetilde{\cL}}|_{(x, - )}$ is surjective. So, there is a dense open set of unfolding parameters $c$ in $(M_{\R})^J$, depending on $\widetilde{\cL}$, such that the tropical sets $P_{\widetilde{\cL}}$ and $p^{-1}_{\widetilde{\cL}}(\del \Gamma_{\cL})$ have codimension one and are transversal in $M_{\R}$. Intersecting these dense open subsets, we can find a value of $c$ such that, for all $\widetilde{\cL}$ belonging to \emph{any fixed, finite subset of} $\LR(\fD_{in, c})$, the tropical sets $P_{\widetilde{\cL}}$ and $p^{-1}_{\widetilde{\cL}}(\del \Gamma_{\cL})$ have codimension one and are transversal in $M_{\R}$. To summarise our discussion, we have
\begin{lem} We can choose an unfolding parameter $c \in (M_{\R})^{J}$ such that, for all $\widetilde{\cL}$ in any fixed finite subset of $\LR(\fD_{in, c})$, there is a dense open subset of the tropical set $P_{\widetilde{L}}$ along which the unique critical point of the potential $W_{\widetilde{\cL}}$ is not contained in the boundary $\del \Gamma_{\cL}$ of the integration cycle $\Gamma_{\cL}$ in the integral formula \eqref{IntegralAlphaUnfolded} for $\alpha_{\widetilde{\cL}}$.
\end{lem} 
Combining this with Corollary \ref{LeadingOrderCor}, and noting that $|\del^2_s W_{\widetilde{\cL},c}| = |\del^2_s W_{\cL}|$, we obtain 
\begin{cor}\label{GenericAsympCor} We can choose an unfolding parameter $c \in (M_{\R})^{J}$ such that, for all $\widetilde{\cL}$ in any fixed finite subset of $\LR(\fD_{in, c})$, there is a dense open subset of $P_{\widetilde{L}}$ along which the leading order term in the $z \to 0$ asymptotic expansion of the $1$-form $\alpha_{\widetilde{\cL}}$ is either irrelevant (i.e. has vanishing transversal integral), or is given by the iterated residue
\begin{equation*}
 2^{(k-1)/2} (\pi(-z))^{-1/2} \frac{e^{W_{\widetilde{\cL},c}(p)/z}}{|\del^2_s W_{\cL}|^{-1/2}} \IR_0(Z_{\widetilde{\cL},c}(\zeta)),
\end{equation*}
where the meromorphic function $Z_{\widetilde{\cL},c}(\zeta)$ is given by 
\begin{equation*}
\frac{r_{\widetilde{\cL}}\Omega_{\widetilde{\cL}}}{\prod^{k-1}_{l = 1} \del_{s_l} W_{\widetilde{\cL},c}(s)} = Z_{\widetilde{\cL},c}(\zeta) {d\zeta_1 \wedge \cdots \wedge d\zeta_{k-1}}.
\end{equation*}
 \end{cor}
 \begin{exm}\label{UnfoldingExm} An unfolding of the integral appearing in Example \ref{ResidueExm2} is given by the sum of $1$-forms $\alpha_{\widetilde{\cL}} + \alpha_{\widetilde{\cL}'}$ where
\begin{align*}
&\alpha_{\widetilde{\cL}}(u) = \left(\frac{1}{\pi (-z)}\right)^{3/2} \int_{(-\infty,0]^2} e^{W_{\widetilde{\cL},c}/z} ds_1 ds_2\,d(2u_1-u_2),\\
&W_{\widetilde{\cL},c} = (u_1 - c_1 -s_2)^2+(u_1 - c_2 -s_1-s_2)^2+(u_2 -s_1-2s_2)^2,\,c_1,\,c_2 \in \R^{3},
\end{align*}
while $\alpha_{\widetilde{\cL}'}$ is obtained from $\alpha_{\widetilde{\cL}}$ by exchanging $c_1, c_2$. 
We compute  
\begin{align*}
& p = (c_1- c_2, \frac{-2c_1+c_2+u_1+u_2}{3}),\, \cW_{\widetilde{\cL}, c}(p) = \frac{1}{3}(2u_1-u_2-c_1-c_2)^2,\\
& \del_{s_1} \cW_{\widetilde{\cL},c} = 2 (c_2+2 s_1+3 s_2-u_1-u_2),\\
&\del_{s_2} \cW_{\widetilde{\cL}, c} = 2 (c_1+c_2+3 s_1 + 6  s_2 -2 u_1 -2 u_2). 
\end{align*}
The same expressions hold for $\alpha_{\widetilde{\cL}'}$ after exchanging $c_1, c_2$. For fixed $c_1 \neq c_2$, the loci $p^{-1}_{\widetilde{\cL}}(\del \Gamma)$, $P_{\widetilde{L}}$ are transversal rays in $M_{\R}$. If $c_1 < c_2$ then the critical point of $W_{\widetilde{\cL}',c}$ is never contained in $\Int(\Gamma)$, and the $1$-form $\alpha_{\widetilde{\cL}'}$ has vanishing transversal integral along $P_{\cL'}$ as $\hbar \to 0$, i.e., it does not contribute to the consistent scattering diagram. On the other hand, on a dense open subset of $P_{\widetilde{L}}$, the critical point of $W_{\widetilde{\cL},c}$ is contained in $\Int(\Gamma)$ and so by Corollary \ref{GenericAsympCor} the leading order term is given by
\begin{align*}
&\left(\frac{1}{2\pi \ii}\right)^2 \int_{C_{\varepsilon}} \frac{ds_1\wedge ds_2}{(2(c_2+2 s_1+3 s_2-u_1-u_2))(2(c_1+c_2+3 s_1 + 6  s_2 -2 u_1 -2 u_2))}\\
&\cdot 2 (\pi(-z))^{-1/2} e^{(2u_1-u_2-c_1-c_2)^2/(3z)} d(2u_1-u_2). 
\end{align*}
Setting
\begin{equation*}
\zeta_1 = 2(c_2+2 s_1+3 s_2-u_1-u_2),\,\zeta_2 = 2(c_1+c_2+3 s_1 + 6  s_2 -2 u_1 -2 u_2) 
\end{equation*}
we compute in this case
\begin{align*}
&Z_{\cL}(\zeta) = \left(\frac{1}{12 \zeta_1 \zeta_2}\right) d(2u_1-u_2),\\
&\alpha_{\cL} \sim 2 (\pi(-z))^{-1/2} (12)^{1/2} e^{(2u_1-u_2-c_1-c_2)^2/(3z)} \IR_{0}(Z_{\cL}(\zeta)).
\end{align*}
Recall that the contribution of $\alpha_{\widetilde{\cL}}$ should be weighted by the factor
\begin{equation*}
g_{\widetilde{\cL}} = (t_1 t_2 g_{\cL})|_{t_1 = t_2 = \frac{1}{2}} = \frac{1}{4} g_{\cL}.
\end{equation*} 
The factor $\frac{1}{4}$ reflects the correct asymptotic behaviour of the original $1$-form $\alpha_{\cL}$ when $\hbar \to 0$ (due to the fact that, before unfolding, we have $p \in \del\Gamma$).
\end{exm} 

\section{JK residues}\label{JKSec}

Here we briefly recall the notion of Jeffrey-Kirwan residues of hyperplane arrangements through their characterisation in terms of flags, following \cite{szenesvergne}, Section 2 (see also \cite{ontani}, Section 4 or \cite{OntaniStoppa}, Section 2 for a detailed exposition), in a very special case, which is sufficient for our purposes. 

Fix a finite set of generators $\mathfrak{G}$ for a lattice $\Lambda$ of finite rank $d$, and suppose it is \emph{projective}, that is, contained in a strict half-space of $\Lambda_{\R}$. We also fix elements $f_1, \ldots, f_d \in \mathfrak{G}$ giving an ordered basis of $\Lambda_{\R}$, and the top form on $\Lambda^{\vee}_{\R}$ defined as $d\mu := f_1 \wedge \cdots \wedge f_d$. Given a subset $S \subset \Lambda$, we will denote by $\mathfrak{B}(S)$ the set of all distinct bases of $\Lambda_{\R}$ consisting of elements of $S$.

Fix a point $x$ in the complexification $\Lambda^{\vee}_{\C} := \Lambda_{\R} \otimes_{\R} \C$. Note that we already considered such a complexification in our discussion of oscillatory integrals, Section \ref{PairingSec}. In that case $d = k-1$ and the lattice $\Lambda$ is generated by the $1$-forms $ds_1,\ldots, ds_{k-1}$. We say that $x \in \Lambda^{\vee}_{\mathbb{C}}$ is \emph{regular}, with respect to the fixed set of generators $\mathfrak{G}$, if it is not contained in the union of hyperplanes
\begin{align*}
\bigcup_{f \in \mathfrak{G}\otimes_\mathbb{R}\mathbb{C}}  V(f) \subset \Lambda^{\vee}_{\mathbb{C}} 
\end{align*}
{($V(f)$ denoting the zero locus of $f$)}. The set of regular points of $\Lambda^{\vee}_{\mathbb{C}}$ is a dense open cone, denoted by $\Lambda^{\vee, reg}_{\mathbb{C}}$.

Similarly, fixing $S \subset \Lambda_{\R}$ a finite set, we say that $\zeta \in \Lambda_{\R}$ is $S$-regular if
		\begin{align*}
			\zeta \notin \bigcup_{\substack{I \subseteq S\\ | I | = d-1}} \operatorname{Span}_\mathbb{R}(I).
		\end{align*}
		In particular, we will say that $\zeta$ is \emph{regular} if it is $\mathfrak{G}$-regular; the locus of such regular points is denoted by $\Lambda^{ reg}_{\R}$. The connected components of the dense open cone $\Lambda^{reg}_{\R}$ are called \emph{$\mathfrak{G}$-chambers of $\Lambda_{\R}$}. 
				
The \emph{Jeffrey-Kirwan residue} we shall consider is a map, depending on the set of generators $\mathfrak{G}$, the volume form $d\mu$ and a chamber $\xi$, defined on the germs at $0 \in \Lambda^{\vee}_{\C}$ of holomorphic functions on the regular locus $\Lambda^{\vee, reg}_{\C}$,
\begin{equation*}
\jk^{\mathfrak{G}}_{\xi}\!: \calH(\Lambda^{\vee, reg}_{\C})_0 \to \C.
\end{equation*} 	 
We will not recall its general definition {(see \cite{brionvergne}, Section 2)}, but only the following characterisation in a special case (see e.g. \cite{OntaniStoppa}, Corollary 2.14, which in turn follows from \cite{szenesvergne}, Theorem 2.6).
\begin{prop}\label{JKBasisProp} Suppose $\mathfrak{G}$ is a basis of $\Lambda_{\R}$. If the chamber $\xi$ is given by
\begin{equation*} 
\xi =\sum_{u \in \mathfrak{G}}\R_{>0}  u
\end{equation*}		
then the map $\operatorname{JK}_{\xi}^{\mathfrak{G}}$ coincides with the composition
\begin{align*}
\calH(\Lambda^{\vee, reg}_{\C})_0 \xrightarrow{\sim} \calH(\C^r)_0 \xrightarrow{\operatorname{IR}_0} \mathbb{C},
\end{align*}
where the first isomorphism is induced by an arbitrary choice of an ordering for the elements of $\mathfrak{G}$ (with corresponding volume form $d\mu$), {while the second map $\operatorname{IR}_0$ is given by the iterated residue of meromorphic functions (see e.g. \cite{OntaniStoppa}, Appendix A).}
	\end{prop}
\subsection{Global residue}\label{GlobalJKSubsec} In the global case, we fix a set $\mathfrak{G}$ of affine functions on $\Lambda^{\vee}_{\R}$, with underlying linear functions induced from a set of generators for $\Lambda$. We say that $x \in \Lambda^{\vee}_{\C}$ is \emph{singular} if it is contained in at least $d$ affine spaces cut out by functions in $\mathfrak{G}$, and write $\mathfrak{G}_x \subset \mathfrak{G}$ for the subset of affine functions vanishing at such $x$. Abusing notation slightly we also write $\mathfrak{G}_x$ for the underlying set of linear functions. Note that a $\mathfrak{G}$-chamber $\xi$ is also automatically a $\mathfrak{G}_{x}$-chamber. The corresponding Jeffrey-Kirwan residue is defined as
\begin{equation*}
\jk^{\mathfrak{G}}_{\xi}(f(s)) = \sum_{x \in \operatorname{Sing}(\mathfrak{G})} \jk^{\mathfrak{G}_x}_{\xi}(f(s + x)),
\end{equation*}    
where the sum is over the set of singular points, and the meromorphic function $f(s)$ is such that $f(s + x)$ is holomorphic in a neighbourhood of the $\mathfrak{G}_x$-regular locus at $s = 0$. 

The definition is extended so that $\jk^{\mathfrak{G}}_{\xi} \equiv 0$ if there are no singular points (i.e. if the linear span of $\mathfrak{G}$ does not generate $\Lambda_{\R}$). 
\section{Application to scattering diagrams}\label{MainThmSec}
The discussion of JK residues in the previous sections generalises at once to functions with values in smooth $1$-forms in the variables $u_1, \ldots, u_r$, as well as in a Lie algebra $\fh$, by linearity. We may now combine Corollary \ref{GenericAsympCor} and Proposition \ref{JKBasisProp} to obtain the following.
\begin{lem} We can choose an unfolding parameter $c \in (M_{\R})^{J}$ such that, for all $\widetilde{\cL}$ in any fixed finite subset of $\LR(\fD_{in, c})$, there is a dense open subset of $P_{\widetilde{L}}$ along which the leading order term in the $z \to 0$ asymptotic expansion of the $1$-form $\alpha_{\widetilde{\cL}}$ is either irrelevant, or is given by the (local) Jeffrey-Kirwan residue
\begin{equation*}
\jk^{\mathfrak{G}}_{\xi}\left(\frac{2^{(k-1)/2} (\pi \hbar)^{-1/2}   e^{-W_{\widetilde{\cL}, c}(p)/\hbar}r_{\widetilde{\cL},c}}{|\del^2_s W_{ \cL }|^{-1/2}\prod^{k-1}_{l = 1} \del_{s_l} W_{\widetilde{\cL},c}(s)}\right),
\end{equation*}
where $\Lambda$ is spanned by $ds_1, \ldots, ds_{k-1}$, $\mathfrak{G} = \{\del_{s_l} W_{\widetilde{\cL},c}\},\,l = 1,\ldots, k-1$ and the chamber $\xi$ is the positive linear span of $\mathfrak{G}$.
\end{lem} 
Note that this also holds (trivially) in the degenerate case when $n_{\cL} = 0$. Applying \cite{leung}, Theorem 3.12 and Proposition 3.14, we then have 
\begin{prop}\label{LocalJKScatterProp} We can choose an unfolding parameter $c \in (M_{\R})^{J}$ such that, for all $\widetilde{L}$ in any fixed finite subset of $\LT(\fD_{in, c})$, the contribution of $\widetilde{L}$ to the completion of $\fD_{in, c}$, which is supported on $P_{\widetilde{L}}$, either vanishes or is given by
\begin{equation*}
\log \Theta_{\widetilde{L}} = -\lim_{\hbar \to 0}\int_{\rho} \frac{1}{|\Aut(\widetilde{L})|}\jk^{\mathfrak{G}}_{\xi}(Z_{\widetilde{L}, c}(s)),
\end{equation*}
for $\rho \subset M_{\R}$ any affine line intersecting positively and generically with $P_{\widetilde{L}}$, where $\mathfrak{G} = \{\del_{s_l} W_{\widetilde{L},c}\},\,l = 1,\ldots, k-1$, the chamber $\xi$ is the positive linear span of $\mathfrak{G}$, and we define
\begin{equation*}
Z_{\widetilde{L}, c}(s) = \frac{2^{(k-1)/2} g_{\widetilde{L}} (\pi \hbar)^{-1/2} e^{-W_{\widetilde{L}, c}(p)/\hbar} r_{\widetilde{L},c}}{|\del^2_s W_{ \cL}|^{-1/2}\prod^{k-1}_{l = 1} \del_{s_l} W_{\widetilde{L},c}(s)}.
\end{equation*}
\end{prop}

\subsection{One-dimensional graded pieces} Under additional assumptions it is possible to emphasise the formal analogy between $Z_{\widetilde{L}, \tau}(s)$ and the rational functions $Z_d(x)$, $Z_{Q, d}(u)$ appearing in Examples \ref{GLSMExm}, \ref{QuiverExm}. 

Suppose that in the graded decomposition 
\begin{equation*}
\fh = \bigoplus_{m \in M^+_{\sigma}} \fh_m
\end{equation*}
each piece is $1$-dimensional, spanned by fixed vectors $h_m$. Note that this is what happens for example in the standard situation when $M$ is endowed with a skew-symmetric bilinear form and $\fh$ is the monoid algebra $\C[M^+_{\sigma}]$ endowed with the corresponding (Kontsevich-Soibelman) Poisson bracket (this includes the particular important case discussed in \cite{gps}, Section 0.1).

In this case, using the definition of $g_{\widetilde{\cL}}$ as an iterated Lie bracket over the cyclically ordered edges of $\widetilde{\cL}$, \cite{leung} Definition 1.14, we can define uniquely scalars $\lambda_l,\,l=1,\ldots, k-1$ attached to $\widetilde{\cL} \in \LR_k$ (depending on the basis $\{h_m\}_{m \in M^+_{\sigma}}$) such that 
\begin{equation*}
g_{\widetilde{\cL}} = \left(\prod^{k-1}_{l = 1} \lambda_l \right) h_{m_{\cL}}. 
\end{equation*}
\begin{cor}\label{KSalgebraCor} In the situation above, the coefficient of $\log \Theta_{\widetilde{L}}$ with respect to the basis vector $h_{m_{\widetilde{L}}}$ is given by
\begin{equation*}
- \lim_{\hbar \to 0} \int_{\rho}\frac{1}{|\Aut(\widetilde{L})|}\jk^{\mathfrak{G}}_{\xi}(\hat{Z}_{\widetilde{L}, c}(s)) ,
\end{equation*}
where
\begin{equation*}
\hat{Z}_{\widetilde{L}, c}(s) =  2^{(k-1)/2} (\pi \hbar)^{-1/2} r_{\widetilde{\cL}} \frac{e^{-W_{\widetilde{L},c}(p)/\hbar}}{|\del^2_s W_{ L} |^{-1/2}}\prod^{k-1}_{l = 1} \left(\frac{\del_{s_l} W_{\widetilde{L},c}(s) + 1}{\del_{s_l} W_{\widetilde{L},c}(s)}\right)^{\lambda_l}.
\end{equation*}
\end{cor}
The proof is an elementary residue computation.
\subsection{Perturbation} Going back to the general situation of Proposition \ref{LocalJKScatterProp}, for each $\widetilde{L}$ in a fixed finite subset of $\fD_{in, c}$, we consider the rational function
\begin{equation*}
Z_{\widetilde{L}, c, \delta_{\widetilde{L}}}(s) := \frac{2^{(k-1)/2} g_{\widetilde{L}} (\pi \hbar)^{-1/2} e^{-W_{\widetilde{L}, c}(p)/\hbar} r_{\widetilde{L},c}}{|\del^2_s W_{L} |^{-1/2}\prod^{k-1}_{l = 1} \del_{s_l} W_{\widetilde{L}, c + \delta_{\widetilde{L}}}(s)},
\end{equation*}
where each $\delta_{\widetilde{L}} \in (M_{\R})^J$ is a perturbation parameter. If $\delta_{\widetilde{L}}$ lies in a sufficiently small neighbourhood of $0 \in (M_{\R})^J$, depending only on the fixed finite subset of $\fD_{in, c}$, then we have
\begin{equation*}
\jk^{\mathfrak{G}}_{\xi}(Z_{\widetilde{L}, c}(s)) = \jk^{\mathfrak{G}_{\delta}}_{\xi_{\delta}}(Z_{\widetilde{L}, c, \delta_{L}}(s)),
\end{equation*} 
where we set $\mathfrak{G}_{\delta} = \{\del_{s_l} W_{\widetilde{L},c+\delta_{L}}\}$ and the chamber $\xi_{\delta}$ is the positive linear span of $\mathfrak{G}_{\delta}$. By choosing the perturbation parameters $\delta_{L}$ sufficiently small and generic, we can ensure that the intersections of affine hyperplanes belonging to any finite subset of $\mathfrak{G}_{\delta}$ are transversal, while preserving the above identity of JK residues for all $\widetilde{L}$ in the fixed finite subset of $\fD_{in, c}$. 

\subsection{Global residue}\label{MainThmSubSec} It is natural to group all labelled trees $\widetilde{L}$ for $\fD_{in, c}$ such that 
\begin{enumerate}
\item[(i)] $P_{\widetilde{L}}$ is a fixed tropical subspace $P \subset M_{\R}$,
\item[(ii)] $g_{\widetilde{L}}$ lies in a fixed $M^+_{\sigma}$-graded component $\fh_m \subset \fh$ (i.e. we have $m_{L} = m$),
\item[(iii)] the transversal integral of $\alpha_{\widetilde{L}}$ along $P$ does not vanish generically. 
\end{enumerate}
We denote this \emph{finite} set by $\cS(P, m)$. Similarly, we define $\cS_k(P, m)$ as the intersection $\cS(P, m) \cap \LT_{k}(\fD_{in, c})$. We also choose $K > 0$ such that $\cS(P, m) = \cS_K(P, m)$.

We would like to express the degree $m$ contribution to the consistent scattering diagram along $P$, namely
\begin{equation*}
- \sum_{\widetilde{L} \in \cS(P, m)} \frac{1}{\Aut(\widetilde{L})}\int_{\rho} \alpha_{\widetilde{L}} g_{\widetilde{L}},
\end{equation*}
as the residue of a suitable meromorphic form.

Consider the rational function
\begin{align}\label{GlobalZ}
\nonumber& Z^{(P, m)}_{c, \delta}(s):= \sum_{k > 0} \sum_{\widetilde{\cL} \in \cS_k(P, m)} \frac{1}{|\Aut(\widetilde{L})|}\frac{Z_{\widetilde{L}, c,\delta_{\widetilde{L}}}(s)}{\prod^K_{j = k} s_j} \\
& = \sum_{k > 0}\sum_{\widetilde{L} \in \cS_k(P, m)}\frac{2^{(k-1)/2}  g_{\widetilde{L}} (\pi \hbar)^{-1/2} e^{-W_{\widetilde{L}, c}(p)/\hbar}r_{\widetilde{L},c}}{|\Aut(\widetilde{L})||\del^2_s W_{L}|^{-1/2}\prod^{k-1}_{l = 1} \del_{s_l} W_{\widetilde{L},c+\delta_{\widetilde{L}}}(s)\prod^K_{j = k} s_j}.
\end{align}
This is holomorphic in the complement of the union of all elements of the set of affine hyperplanes (all contained in the same affine space, by using the natural inclusions) given by
\begin{equation*}
\mathfrak{G}^{(P, m)}_{\delta}:= \{\del_{s_l} W_{\widetilde{L},c + \delta_{\widetilde{L}}}\!:\widetilde{L} \in \cS_k(P, m),\, l = 1,\ldots, k-1\} \bigcup \{s_j = 0,\,j = 1,\cdots, K\}.
\end{equation*}
For all small and sufficiently general perturbation parameters $\delta_{\widetilde{L}}$, all intersections of all subsets of $\mathfrak{G}^{(P, m)}_{\delta}$ are transverse. In particular, if we let the chamber $\xi^{(P, m)}_{\delta}$ be given by the positive cone spanned by elements of $\mathfrak{G}^{(P, m)}_{\delta}$, then by Proposition \ref{JKBasisProp} and our discussion in Section \ref{GlobalJKSubsec}, we have
\begin{align*}
&\jk^{\mathfrak{G}^{(P, m)}_{\delta}}_{\xi^{(P, m)}_{\delta}}(Z^{(P,m)}_{c,\delta}) = \sum_{x \in \operatorname{Sing}(\mathfrak{G}^{(P, m)}_{ \delta })} \jk^{\mathfrak{G}^{(P, m)}_{ \delta, x}}_{\xi^{(P, m)}_{\delta}}(Z^{(P,m)}_{c, \delta}(s +x))\\
&= \sum_{x \in \operatorname{Sing}(\mathfrak{G}^{(P, m)}_{\delta})} \IR_0(Z^{(P,m)}_{c, \delta}(s+x)).
\end{align*}
By transversality, the sum is actually taken over intersections $x = \bigcap^{k-1}_{i = 1} H_{i}$ of distinct elements $H_1, \ldots, H_{k-1}$ of $\mathfrak{G}^{(P, m)}_{\delta}$.

Let $Z_{\widetilde{L}, c, \delta_{\widetilde{L}}}(s)$ be a summand of $Z^{(P, m)}_{c, \delta}(s)$ in the decomposition \eqref{GlobalZ}, and fix a singular point $x = \bigcap^{k-1}_{i = 1} H_{i}$ as above. Then, if the sets $\{H_{i},\,i = 1,\ldots,k-1\}$ and $\{\del_{s_l} W_{\widetilde{L},c+\delta_{\widetilde{L}}},\, l = 1,\ldots, k-1\}$ do not coincide, there is at least one $H_j$ such that the meromorphic function $Z^{(P, m)}_{k, c,\delta}(s + x)$ is regular along $H_j$ in a neighbourhood of $s = 0$, and so the corresponding iterated residue vanishes. Therefore, we have 
\begin{equation*}
\IR_0(Z^{(P,m)}_{c}(s+x)) = \IR_0(Z_{\widetilde{L}(x), c}(s+x)),
\end{equation*}
where $\widetilde{L}(x)$ is the unique labelled tree such that the sets $\{H_{i},\,i = 1,\ldots,k-1\}$ and $\{\del_{s_l} W_{\widetilde{L},c+\delta_{\widetilde{L}}},\, l = 1,\ldots, k-1\}$ coincide. In turn, we have
\begin{equation*}
\IR_0(Z_{\widetilde{L}(x), c, \delta_{\widetilde{L}(x)}}(s+x)) = \jk^{\mathfrak{G}^{(P, m)}_{ \delta }}_{\xi^{(P, m)}_{ \delta }}(Z_{\widetilde{L}, c, \delta_{\widetilde{L}}}(s)),
\end{equation*}
and by using Proposition \ref{LocalJKScatterProp} we obtain
\begin{thm}\label{MainThm} For suitable unfolding parameters $c \in (M_{\R})^{J}$, the contribution of the set of all labelled trees $\widetilde{L} \in \cS(P, m)$ to the completion of $\fD_{in, c}$, which is supported on $P$, is given by
\begin{equation*}
 -\lim_{\hbar \to 0}\int_{\rho} \jk^{\mathfrak{G}^{(P, m)}_{\delta}}_{\xi^{(P, m)}_{\delta}}(Z^{(P,m)}_{ c, \delta}),
\end{equation*}
where the perturbation parameters $\delta = \{\delta_{L}\}$ are sufficiently small and generic, and $\rho \subset M_{\R}$ is any affine line intersecting positively and generically with $P$. 
\end{thm} 
Note that, by Corollary \ref{KSalgebraCor}, in the special case when each graded piece $\fh_m$ is $1$-dimensional, with a distinguished basis vector $h_m$, we may write 
\begin{align*} 
& Z^{(P, m)}_{c, \delta}(s)=\\ 
& \sum_{k > 0}\sum_{\widetilde{L} \in \cS_k(P, m)}\frac{2^{(k-1)/2} (\pi \hbar)^{-1/2} e^{-W_{\widetilde{L}, c}(p)/\hbar}r_{\widetilde{L},c}}{|\Aut(\widetilde{L})||\del^2_s W_{ L}|^{-1/2}\prod^K_{j = k} s_j} \prod^{k-1}_{l = 1} \left(\frac{\del_{s_l} W_{\widetilde{L},c+\delta_{\widetilde{L}}}(s) + 1}{\del_{s_l} W_{\widetilde{L},c+\delta_{\widetilde{L}}}(s)}\right)^{\lambda_{\widetilde{L}, l}} h_{m_{\widetilde{L}}}
\end{align*}
for unique scalars $\{\lambda_{\widetilde{L}, l}\}$.
{
\begin{exm}\label{QuiverExm1} Let us consider the simplest application to quivers. Following Example \ref{KroneckerExm}, we consider the contribution of framed representations of $Q_{\kappa}$ of dimension vector $(1, 1)$ to the consistent scattering diagram $\fD$. According to Theorem \ref{LMYThm}, the coefficient $c_{(1,1)}$, that is, the degree $(1,1)$ contribution to the (infinitesimal generator of the) wall-crossing automorphism $\log \Theta_P$ in the wall
\begin{equation*}
\{m := (1,1),\,n = \{ -, m\},\,P = \R_{>0} m,\,\Theta_P \}
\end{equation*}
is given by
\begin{equation*}
c_{(1,1)} = -\lim_{\hbar \to 0}\int_{\rho}\alpha_{L} g_{L},
\end{equation*}
where $L$ is the unique tree with roots labelled by $x = x_{m_1}$, $y = x_{m_2}$ (in the notation of Example \ref{KroneckerExm}). By Example \ref{ResidueExm1}, we can compute 
\begin{align*}
& -\lim_{\hbar \to 0}\int_{\rho}\alpha_{L} g_{L} = [x, y]\int_{\rho} (-\opH(\delta \wedge \eta))\\
& = \kappa\,x y \lim_{z \to 0}\int_{\rho} 2\sqrt{2}(\pi(-z))^{-1/2} \jk^{\mathfrak{G}^{(P, m)}_{ 0 }}_{\xi^{(P, m)}_{ 0 }}\left(\frac{1}{-2(u_1+u_2 - 2s)}\right) e^{\frac{(u_1-u_2)^2}{2 z}} d(u_1 - u_2)\\
& = \kappa\,x y. 
\end{align*}
This is consistent with the fact that $\M^{(1, 0), B}_{\kappa}(1, 1) \cong \M^{(0, 1), F}_{\kappa}(1, 1) \cong \PP^{\kappa - 1}$. Note that in this case, as explained in Example \ref{ResidueExm1}, unfolding is not required. 
\end{exm}
}
{
\begin{exm}\label{QuiverExm2} For a slightly more complicated example we consider the contribution of framed representations of $Q_{\kappa}$ of dimension vector $(1, 2)$ to the consistent scattering diagram $\fD$. According to Example \ref{KroneckerExm} and Theorem \ref{LMYThm}, the coefficient $c_{(1,2)}$, that is, the degree $(1,2)$ contribution to the (infinitesimal generator of the) wall-crossing automorphism $\log \Theta_P$ in the wall
\begin{equation*}
\{m := (1,2),\,n = \{ -, m\},\,P = \R_{>0} m,\,\Theta_P \}
\end{equation*}
is given by
\begin{equation*}
c_{(1, 2)} = -\lim_{\hbar \to 0}\int_{\rho}\alpha_{L} g_{L} -\lim_{\hbar \to 0}\int_{\rho}\alpha_{L'} g_{L'},
\end{equation*}
where $L$, $L'$ are the unique trees with roots labelled by $\{x, y, y\}$, respectively $\{x, \frac{y^2}{4}\}$, such that $g_{L} \neq 0$. By Examples \ref{ResidueExm2} and \ref{UnfoldingExm}, we can compute
\begin{align*}
&- \int_{\rho}\alpha_{L} g_{L} = [[x, y],y] \int_{\rho} (-\opH(-\opH(\delta \wedge \eta) \wedge \eta))\\
& \sim \frac{\kappa^2}{4} x y^2\left(\frac{1}{2\pi \ii}\right)^2  2 (\pi(-z))^{-1/2} \jk^{\mathfrak{G}^{(P, m)}_{ 0 }}_{\xi^{(P, m)}_{ 0 }}(Z(s)) \int_{\rho} e^{(2u_1-u_2-c_1-c_2)^2/(3z)} d(2u_1-u_2),
\end{align*} 
where
\begin{equation*}
Z(s) = (2(c_2+2 s_1+3 s_2-u_1-u_2))^{-1}(2(c_1+c_2+3 s_1 + 6  s_2 -2 u_1 -2 u_2))^{-1}.
\end{equation*}
Here $c_1$, $c_2$ are unfolding parameters satisfying $c_1 < c_2$, and we are choosing the perturbation parameter $\delta = 0$. 

On the other hand, with computations entirely similar to Example \ref{ResidueExm1}, but keeping track of the unfolding parameters, we find  
\begin{align*}
& -\int_{\rho}\alpha_{L'} g_{L'} =  -[x, \frac{y^2}{4}] \int_{\rho}\sum^2_{i = 1}\frac{1}{2}\left(\frac{1}{\pi (-z)}\int^0_{-\infty} e^{((u_1- c_i -s_1)^2+(u_2-2s_1)^2)/z} ds \right) d(2 u_1-u_2)\\
& \sim -\frac{\kappa}{2} x y^2 2^{1/2}(\pi(-z))^{-1/2} (10)^{1/2} \sum^2_{i=1}\frac{1}{2}\jk^{\mathfrak{G}^{(P, m)}_{ 0 }}_{\xi^{(P, m)}_{ 0 }}(Z'_i(s)) \int_{\rho} e^{\frac{(2u_1-u_2-2c_i)^2}{5 z}}d(2 u_1-u_2), 
\end{align*}
where 
\begin{equation*}
Z'_i(s) = \frac{1}{-2 (-c_i - 5 s_1 + u_1 + 2 u_2) s_2}.
\end{equation*}
Evaluating the Gaussian integrals yields the JK residue representation  
\begin{equation*}
c_{(1, 2)} = x y^2\jk^{\mathfrak{G}^{(P, m)}_{ 0 }}_{\xi^{(P, m)}_{ 0 }}\left(\frac{\kappa^2}{2} Z(s) - \frac{\kappa}{4} Z'_1(s) - \frac{\kappa}{4} Z'_2(s)\right),   
\end{equation*}
which is valid for generic perturbation parameters $c_1$, $c_2$ satisfying $c_1 < c_2$. Note that computing the residues we find $c_{(1, 2)} = (\frac{\kappa^2}{2} - \frac{\kappa}{2}) x y^2$. In particular, for $\kappa =1$ we have $c_{(1, 2)} = 0$, which must hold since $\M^{(1, 0), B}_{1}(1, 2) = \M^{(0, 1), F}_{1}(1, 2) = \emptyset$; while for $\kappa = 2$, we have $c_{(1, 2)} = x y^2$, which must be the case since $\M^{(1, 0), B}_{2}(1, 2) \cong \{pt\}$, $ \M^{(0, 1), F}_{2}(1, 2) \cong \PP^1$ (see \cite{gp}, Section 2.6 for both statements).
\end{exm}
}
\section{Theta functions}\label{ThetaSec}

The Maurer-Cartan scattering approach can be extended to construct the \emph{theta functions} $\vartheta_{m}$ of the consistent completion $\fD$, with values in a tropical representation $A$ of $\fh$ in the sense of \cite{leung}, Section 1.1.1.

Following \cite{leung}, Section 1.3, the class of labelled ribbon $k$-trees $\LR_k$ (with respect to $\fD_{in}$) is extended to the class of \emph{marked} labelled ribbon $k$-trees $\mr_k$, containing one distinguished edge $\check{e}$. The directed path starting from $\check{e}$ is the \emph{core} $\fc$ of $\cJ \in \mr_k$. Excising $\fc$ from $\cJ$ produces labelled ribbon trees $\cL_1, \ldots, \cL_l$ for some $l$ (the length of $\fc$), with underlying labelled trees $L_1, \ldots, L_l$. The distinguished Maurer-Cartan solution \eqref{KuranishiMCelement} is replaced by the algebra element
\begin{equation*}
\theta_m = \sum_{k > 0} \frac{1}{2^{k-1}} \sum_{\cJ \in \mr_k}\opL_{k, \cJ}(\Pi + z^{\varphi(m)}, \ldots, \Pi + z^{\varphi(m)}), 
\end{equation*}  
where the operator $\opL_{k, \cJ}$ extends the previous definition by placing $\Pi$ at unmarked edges and $z^{\varphi(m)}$ at $\check{e}$. According to \cite{leung} (equation after (3.6)) we have in terms of marked trees $J \in \mt_k$ (forgetting the ribbon structure)
\begin{equation*}
\theta_m = \sum_{k > 0} \sum_{J \in \mt_k} \frac{1}{|\Aut(J)|} \alpha_J a_J,
\end{equation*}
for certain algebra elements $a_J \in A$ and smooth \emph{functions} $\alpha_J$ (i.e. $0$-forms rather than $1$-forms as in the case of $\alpha_L$).

The relation to the theta functions of $\fD$ is given by
\begin{equation*}
\lim_{\hbar \to 0} \theta_m = \vartheta_m
\end{equation*}
(see \cite{leung} Theorem 3.20). 

Following the proof of \cite{leung}, Lemma 3.19, for each $J \in \mt_k$, one defines a vector of linear flows $\tau_{\fc} = (\tau_{\fc, 1}, \ldots, \tau_{\fc, l})$ (where $l$ is the length of the core $\fc$ of $J$), such that 
\begin{equation*}
\alpha_J(x) = (-1)^l \int_{(-\infty, 0]^l \times\{x\}} \tau^*_{\fc}(\alpha_{L_1} \wedge\cdots\wedge\alpha_{L_l}). 
\end{equation*}  
As in \eqref{1formAlpha} we can write, for $i = 1, \dots, l$,
\begin{align*}  
& \alpha_{L_i}(x)=\\
& \left(\frac{1}{\pi \hbar}\right)^{k_{L_i}/2}\int_{(-\infty, 0]^{L^{[1]}_i}\times \{x\}} e^{-\big(\sum^{k_{L_i}}_{j=1} (\tau^{\fe_{L_i,j}}_{L_i})^* \eta^2_{i_{e_{L_i,j}}}\big)/\hbar} \vec{\tau}^*_{L_i}\big( d\eta_{e_{L_i,1}} \wedge \cdots \wedge d\eta_{e_{L_i,k}} \big).
\end{align*} 
Therefore, we also have
\begin{align*}
\alpha_J = (-1)^l \left(\frac{1}{\pi \hbar}\right)^{k_J/2} \int_{\Gamma_{J}} e^{-W_J(s)/\hbar} r_J \Omega_J(s)
\end{align*}
where
\begin{align*}
& k_J = \sum^l_{i = 1} k_{L_i},\, \Gamma_J = (-\infty, 0]^l \times (-\infty, 0]^{L^{[1]}_1} \times \cdots \times (-\infty, 0]^{L^{[1]}_l},\\
& W_J(s) = \sum^l_{i=1}\tau^*_{\fc, i}\left(\sum^{k_{L_i}}_{j=1} (\tau^{\fe_{L_i,j}}_{L_i})^* \eta^2_{i_{e_{L_i,j}}}\right),\\
&\Omega_J(s)= ds_1\wedge\cdots\wedge ds_l \wedge \bigwedge^l_{i=1}\left(\bigwedge^{k_{L_i} -1 }_{j = 1}  d s_{L_i, j}\right),\\
&r_J = \frac{\tau^*_{\fc}\big(\bigwedge^l_{i=1}\vec{\tau}^*_{L_i}\big( d\eta_{e_{L_i,1}} \wedge \cdots \wedge d\eta_{e_{L_i,k}} \big)\big)}{\Omega_J}.
\end{align*} 

Fix an unfolded initial diagram $\fD_{in, c}$. We denote by $\cS(Q, m, \mu)$ the set of $z^{\varphi(m)}$-marked trees $\widetilde{J}$ \emph{for the unfolded scattering diagram} $\fD_{in, c}$ such that the corresponding algebra element $a_{\widetilde{J}}$ has fixed degree $\mu \in M^+_{\sigma}$, and which give a \emph{nonvanishing contribution to the theta function $\widetilde{\vartheta}_m$ for the consistent completion $\fD_{c}$ of $\fD_{in, c}$}, evaluated at the point $Q \in M_{\R} \setminus \fD_{c}$. Note that the latter condition implies in particular that, letting $p_{\widetilde{J}}$ be the unique critical point of the quadratic potential $W_{\widetilde{J}}(s)$, we have
\begin{equation*} 
W_{\widetilde{J}}(p_{\widetilde{J}}) = 0.
\end{equation*}

For each $\widetilde{J} \in \cS(Q, m, \mu)$, we write
\begin{equation*}
\operatorname{Ind}(\widetilde{J}) = \{1, \ldots, l\}\times\{1, \ldots, k_{L_i}-1\}\times \cdots\times\{1, \ldots, k_{L_l}-1\}, 
\end{equation*} 
where $l$ denotes the length of the core $\fc$ of $\widetilde{J}$ as usual. We think of all these index sets as contained in a sufficiently large finite set $\operatorname{Ind}(\mu)$.

We introduce rational functions, with values in $A$, 
\begin{equation*}
Z_{\widetilde{J}, c, \delta_{\widetilde{J}}}(s) := \frac{2^{(k_J -1)/2} r_{\widetilde{J},c} a_{\widetilde{J}}}{|\del^2_s W_{J}|^{-1/2}\prod_{q \in \operatorname{Ind}(\widetilde{J})} \del_{s_q} W_{\widetilde{J}, c + \delta_{\widetilde{J}}}(s)},
\end{equation*}
where $\delta_{\widetilde{J}}$ denote sufficiently small and generic perturbation parameters.

Similarly, we define rational functions
\begin{align*} 
\nonumber& Z^{(Q, m, \mu)}_{c, \delta}(s) := \sum_{\widetilde{J} \in \cS(Q, m, \mu)} \frac{1}{|\Aut(\widetilde{J})|}\frac{Z_{\widetilde{J}, c, \delta_{\widetilde{J}}}(s)}{\prod_{q \in \operatorname{Ind}(\mu)\setminus \operatorname{Ind}(\widetilde{J})} s_q} \\
& = \sum_{\widetilde{J} \in \cS(Q, m, \mu)} \frac{1}{|\Aut(\widetilde{J})|} \frac{2^{(k_J -1)/2} r_{\widetilde{J},c} a_{\widetilde{J}}}{|\del^2_s W_{J}|^{-1/2}\prod_{q \in \operatorname{Ind}(\widetilde{J})} \del_{s_q} W_{\widetilde{J}, c + \delta_{\widetilde{J}}}(s)\prod_{q \in \operatorname{Ind}(\mu) \setminus \operatorname{Ind}(\widetilde{J})} s_q}.
\end{align*}

Now the arguments of Sections \ref{UnfoldingSec} and \ref{MainThmSec} go through with minor changes, and yield the following.
\begin{prop}\label{ThetaProp} For fixed $m$ and generic $c$, $\delta$ and $Q \in M_{\R} \setminus \fD_{c}$, we have
\begin{equation*}
\widetilde{\vartheta}_m(Q) = \sum_{\mu \in M^+_{\sigma}} \jk^{\mathfrak{G}^{\mu}_{\delta}}_{\xi^{\mu}_{\delta}}(Z^{(Q, m, \mu)}_{c, \delta}(s)),
\end{equation*}
where each set of affine hyperplanes $\mathfrak{G}^{\mu}_{\delta}$ is given by
\begin{equation*}
\mathfrak{G}^{\mu}_{\delta}:= \{\del_{s_q} W_{\widetilde{J}, c + \delta_{\widetilde{J}}}\!:\widetilde{J} \in \cS(Q, m, \mu),\, q \in \operatorname{Ind}(\widetilde{J})\} \bigcup \{s_q = 0,\,q \in \operatorname{Ind}(\mu) \setminus \operatorname{Ind}(\widetilde{J})\} 
\end{equation*}
and the chamber $\xi^{\mu}_{\delta}$ is the positive linear span of $\mathfrak{G}^{\mu}_{\delta}$.
\end{prop}
In the special case when each graded piece $\fh_m$ is $1$-dimensional, with a distinguished basis vector $h_m$, we may write
\begin{align*} 
& Z^{(Q, m, \mu)}_{c, \delta}(s)= \sum_{\widetilde{J} \in \cS(Q, m, \mu)} \frac{1}{|\Aut(\widetilde{J})|} \frac{2^{(k_J -1)/2} r_{\widetilde{J},c} \hat{a}_{\widetilde{J}}}{|\del^2_s W_{J}|^{-1/2}\prod_{q \in \operatorname{Ind}(\mu) \setminus \operatorname{Ind}(\widetilde{J})} s_q}\\&\quad\quad\quad\quad\quad\quad\,\,\,\,\prod_{q \in \operatorname{Ind}(\widetilde{J})}\left(\frac{\del_{s_q} W_{\widetilde{J}, c + \delta_{\widetilde{J}}}(s)+1}{ \del_{s_q} W_{\widetilde{J}, c + \delta_{\widetilde{J}}}(s)}\right)^{\lambda_q},
\end{align*}
for suitable scalars $\lambda_q$, such that
\begin{align*}
& g_{\tilde{L}_i} = \left(\prod^{k_i}_{j= 1} \lambda_{L_i, j} \right) h_{m_{\tilde{L}_i}}\,i = 1,\ldots, l;\,\,\hat{a}_{\widetilde{J}} = h_{m_{\tilde{L}_l}} \cdot \,(\cdots)\, \cdot h_{m_{\tilde{L}_1}} z^{\varphi{(m)}}.
\end{align*}

\end{document}